\documentclass{amsart}
\usepackage{amsmath}
\usepackage{amsthm} 
\usepackage{amssymb} 
\usepackage{amscd} 
\usepackage{graphicx} 
\usepackage{epsfig}
\usepackage[dvips]{color}

\renewcommand{\labelenumi}{ $(\arabic{enumi})$ }

\newtheorem{theorem}{Theorem}[section]
\newtheorem{lemma}[theorem]{Lemma}

\newtheorem{proposition}[theorem]{Proposition}
\newtheorem{corollary}[theorem]{Corollary}
\newtheorem{claim}[theorem]{Claim}

\theoremstyle{definition}
\newtheorem{definition}[theorem]{Definition}
\newtheorem{example}[theorem]{Example}

\theoremstyle{remark}
\newtheorem{remark}[theorem]{Remark}
\newtheorem{question}[theorem]{Question}
\newtheorem{conjecture}[theorem]{Conjecture}

\numberwithin{equation}{section}
\numberwithin{figure}{section}
\numberwithin{table}{section}

\newcommand{\QED}[1]{\hspace*{\fill} $\square$(#1)\newline}

\begin{document}

\title
{Networking Seifert Surgeries on Knots III}

\author{Arnaud Deruelle,\ 
Katura Miyazaki
and Kimihiko Motegi \\
}

\date{}
\dedicatory{Dedicated to Sadayoshi Kojima on the occasion of his 60th birthday}

\begin{abstract}
How do Seifert surgeries on hyperbolic knots arise from those on torus knots? 
We approach this question from a networking viewpoint introduced in \cite{DMM1}. 
The Seifert Surgery Network is
a $1$--dimensional complex
whose vertices correspond to Seifert surgeries;
two vertices are connected by an edge if
one Seifert surgery is obtained from the other by
a single twist along a trivial knot called a seiferter
or along an annulus cobounded by seiferters.
Successive twists along a ``hyperbolic seiferter"
or a ``hyperbolic annular pair" produce
infinitely many Seifert surgeries on hyperbolic knots.
In this paper, we investigate Seifert surgeries on torus knots
which have hyperbolic seiferters or hyperbolic annular pairs,
and obtain results suggesting that
such surgeries are restricted.

\end{abstract}

{
\renewcommand{\thefootnote}{}
}
 
\maketitle

\section{Introduction}
\label{section:Introduction}
\textit{How do Seifert surgeries on hyperbolic knots arise from those on torus knots?}
In \cite{DMM1} we formulate this question from a viewpoint
of the Seifert Surgery Network. 
Let us recall some basic notions given in \cite{DMM1} and 
an example illustrating our idea.  
A pair $(K, m)$ of a knot $K$ in $S^3$ and an integer $m$ is 
a \textit{Seifert surgery} if the result $K(m)$ of $m$--Dehn surgery on $K$ has a Seifert fibration;
we allow the fibration to be degenerate,
i.e.\ it contains an exceptional fiber of index 0 as a degenerate fiber.
It is shown in \cite[Proposition~2.8]{DMM1} that
if $K(m)$ admits a degenerate Seifert fibration,
then it is either a lens space or
a connected sum of two lens spaces. 
In the latter case, Greene \cite{Greene} recently shows that 
$K$ is a torus knot or a cable of a torus knot. 

\begin{definition}[\textbf{seiferter}]
Let $(K, m)$ be a Seifert surgery.
A knot $c$ in $S^3 -N(K)$ is called a \textit{seiferter} for $(K, m)$
if $c$ satisfies (1) and (2) below.
\begin{enumerate}
\item $c$ is a trivial knot in $S^3$.
\item $c$ becomes a fiber in
a Seifert fibration of $K(m)$.
\end{enumerate}
\end{definition}

We also consider pairs of seiferters.

\begin{definition}[\textbf{annular pair of seiferters}]
\label{annular pair}
Let $c_1,c_2$ be seiferters for a Seifert surgery $(K, m)$.
We call $\{ c_1, c_2 \}$ a \textit{pair of seiferters} if 
$c_1$ and $c_2$ simultaneously become fibers in a Seifert fibration
of $K(m)$. 
A pair of seiferters $\{c_1, c_2\}$ is called a \textit{Hopf pair}
if $c_1 \cup c_2$ is a Hopf link in $S^3$.
A pair of seiferters $\{c_1, c_2\}$ is called an \textit{annular pair of seiferters} 
(or \textit{annular pair} for short) if 
$c_1$ and $c_2$ cobound an annulus in $S^3$. 
\end{definition}

For a Seifert surgery $(K, m)$ with a seiferter $c$,
let $K_p$ and $m_p$ be the images of $K$ and $m$
under $p$--twist along $c$, respectively.
Then, $(K_p, m_p)$ remains a Seifert surgery for any integer $p$,
and (the image of) $c$ is also a seiferter
for $(K_p, m_p)$ (\cite[Proposition~2.6]{DMM1}).
Similarly, if $(K, m)$ has an annular pair $\{ c_1, c_2 \}$, 
then under twisting along the annulus cobounded by $c_1, c_2$, 
we obtain a new Seifert surgery for which (the image of) $\{c_1, c_2\}$ remains an annular pair
 (\cite[Proposition~2.33(1)]{DMM1}). 
We call a twist along an annulus cobounded by $c_1 \cup c_2$
a \textit{twist along an annular pair of seiferters}.
We say that
a seiferter $c$ (resp.\ an annular pair $\{c_1, c_2\}$)
for a Seifert surgery $(K, m)$ is \textit{hyperbolic}
if $S^3 - K \cup c$ (resp.\ $S^3 - K \cup c_1 \cup c_2$)
admits a complete, hyperbolic metric with finite volume. 

\begin{remark}
\label{rem:irrelevant}
Suppose that a seiferter $c$ for $(K, m)$ bounds a disk in $S^3 -K$.
Since no twist along $c$ changes $(K, m)$,
we call $c$ \textit{irrelevant}.
\textit{We do not regard an irrelevant seiferter as a seiferter}.
However,
for pairs of seiferters $\{c_1, c_2\}$ we allow $c_i$ to be an irrelevant seiferter.
Let $\{c_1, c_2\}$ be an annular pair for $(K, m)$.
If either $c_1$ and $c_2$ cobound an annulus disjoint from $K$ or
there is a 2--sphere in $S^3$ separating $c_i$ and $c_j\cup K$,
then twists along $\{c_1, c_2\}$ do not change $(K, m)$ or
have the same effect on $K$ as twists along $c_j$.
We thus call such an annular pair \textit{irrelevant},
and exclude it from annular pairs of seiferters.
Note that if $S^3 -K\cup c_1\cup c_2$ is hyperbolic,
then $\{c_1, c_2\}$ is not irrelevant.
\end{remark}

Regard each Seifert surgery as a vertex,
and connect two vertices by an edge if
one is obtained from the other by a single twist along a seiferter
or an annular pair of seiferters. 
We then obtain a $1$--dimensional complex, called the \textit{Seifert Surgery Network}. 

Let us take a look at seiferters for Seifert surgeries on torus knots $T_{p, q}$. 
Throughout this paper we assume, 
without loss of generality, 
that $|p| > q \ge 1$, 
and denote a trivial knot $T_{p, 1}$ by $O$. 

\begin{figure}[htbp]
\begin{center}
\includegraphics[width=0.29\linewidth]{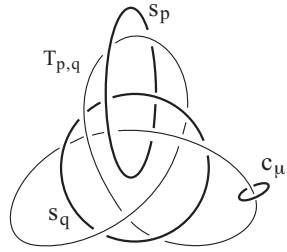}
\caption{Basic Seiferters}
\label{basicseiferters}
\end{center}
\end{figure}

\begin{example}[\textbf{the subcomplex $\mathcal{T}$}]
\label{basic seiferters}
Since the exterior of a torus knot $T_{p, q}$ is a Seifert fiber space, 
$(T_{p, q}, m)$ is a Seifert surgery for any integer $m$. 
Let $s_p$, $s_q$ be exceptional fibers in 
the Seifert fibration of the exterior of $T_{p,q}$
with indices $|p|$, $q$, respectively,
and $c_{\mu}$ a meridian of $T_{p,q}$;
see Figure~\ref{basicseiferters}.
Then $s_p$ and $s_q$
remain exceptional fibers in $T_{p, q}(m)$. 
Note that $c_{\mu}$ is isotopic in $T_{p,q}(m)$
to the core of the filled solid torus,
which is a fiber of index $|pq -m|$
and in particular a degenerate fiber in $T_{p, q}(pq)$. 
Hence, the trivial knots
$s_p$, $s_q$, $c_{\mu}$ are seiferters for $(T_{p,q}, m)$
for any integer $m$, 
and called \textit{basic seiferters} for $(T_{p,q}, m)$. 
We denote by $\mathcal{T}$
the subcomplex such that
its vertices are Seifert surgeries on torus knots and
its edges correspond to basic seiferters.  
\end{example}

The following example motivates us to consider the Seifert Surgery Network. 

\begin{figure}[htbp]
\begin{center}
\includegraphics[width=0.95\linewidth]{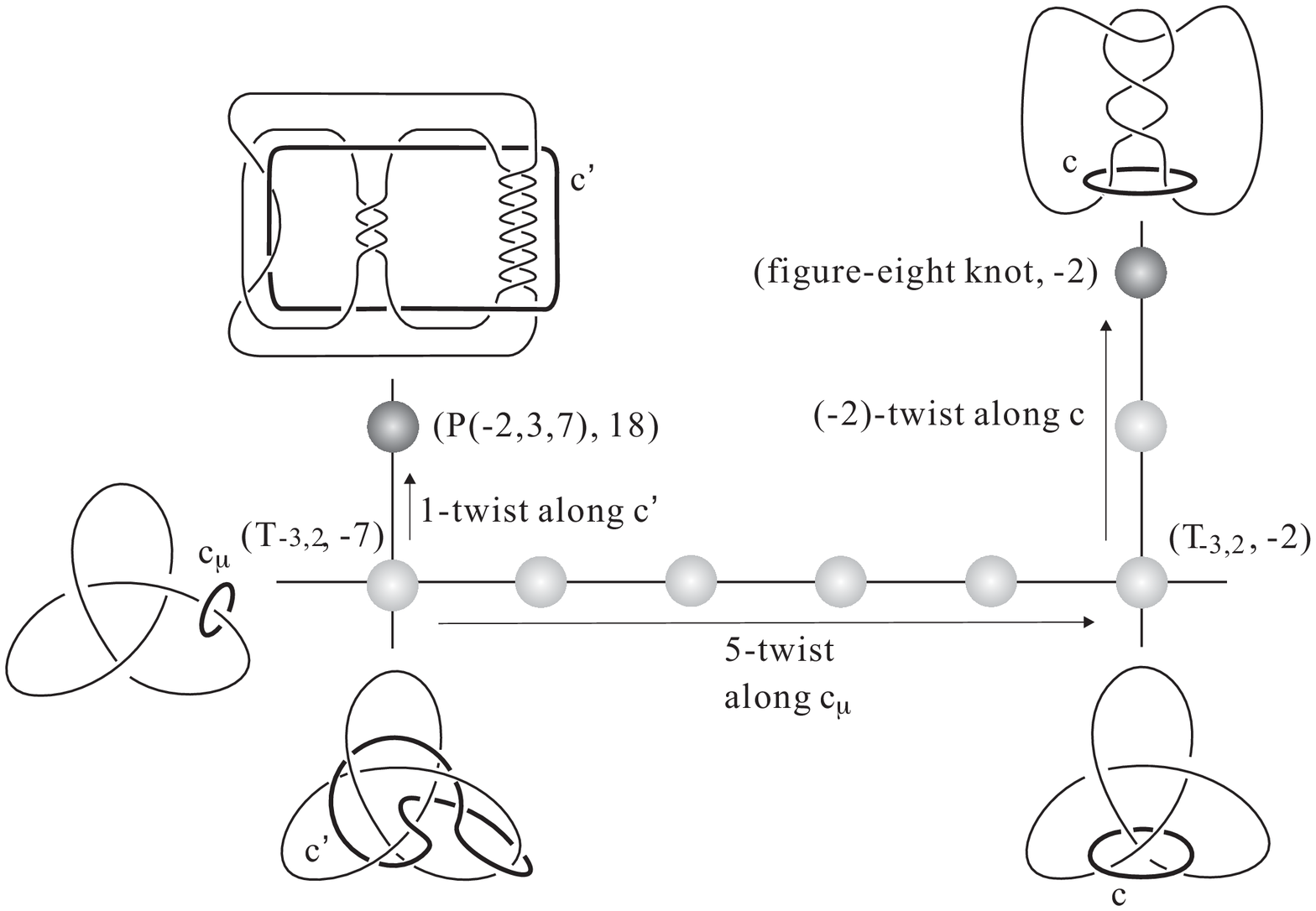}
\caption{Seifert Surgery Network}
\label{motivationalExample}
\end{center}
\end{figure}

\begin{example}
\label{motivation}
\begin{enumerate}
\item
The meridian $c_{\mu}$ for $T_{-3, 2}$ is a seiferter for
all $(T_{-3, 2}, m)$ $(m\in \mathbb{Z})$.
Twisting along $c_{\mu}$ yields the horizontal line in Figure~\ref{motivationalExample}, 
which consists of all the integral 
Seifert surgeries on $T_{-3, 2}$.  

\item
The trivial knot $c \subset S^3 -T_{-3, 2}$
in Figure~\ref{motivationalExample}
is a seiferter for the Seifert surgery $(T_{-3, 2}, -2)$
(Section~\ref{m-move}, Figure~\ref{minus2move}).  
A $(-2)$--twist of $T_{-3, 2}$ along $c$ yields 
the figure--eight knot. 
Since the linking number between $c$ and $T_{-3, 2}$ is zero, 
the surgery slope $-2$ does not change under the twisting. 
Thus we obtain the right vertical line in Figure~\ref{motivationalExample}.

\item
The trivial knot $c' \subset S^3 -T_{-3, 2}$
in Figure~\ref{motivationalExample} is a seiferter for
$(T_{-3, 2}, -7)$ (\cite[Example~2.21(2)]{DMM1}).
A $1$--twist of $T_{-3, 2}$ along $c'$ yields
the $(-2, 3, 7)$ pretzel knot $P(-2, 3, 7)$ (\cite{DMMtrefoil}).
Since the linking number between $c'$ and $T_{-3, 2}$ is $5$,
the surgery slope changes from $-7$ to $-7 +5^2 =18$.  
We thus obtain the lens surgery $(P(-2, 3, 7), 18)$
first found by Fintushel and Stern \cite{FS}. 
This gives the left vertical line
in Figure~\ref{motivationalExample}. 
\end{enumerate}
\end{example}

A path from $(K, m)$ to $(K', m') \in \mathcal{T}$
in the network shows that
the Seifert surgery $(K, m)$ is obtained from
the $m'$--surgery on the torus knot $K'$ by a sequence of twists
along seiferters or annular pairs.
For example, 
vertical paths in Figure~\ref{motivationalExample} 
from $($figure-eight knot$, -2)$ and $(P(-2,3,7), 18)$
to vertices in $\mathcal{T}$
explain how these surgeries arise from 
surgeries on a trefoil knot. 
In \cite{DMM2, DMM1, DEMM},
we find paths from various known Seifert surgeries
to vertices in $\mathcal{T}$;
the list includes Seifert surgeries on graph knots,
Berge's lens surgeries \cite{Berge2}, 
and Seifert surgeries constructed by using Montesinos trick \cite{EM1, EM2}. 

In the present paper, 
we explore Seifert surgeries on torus knots which have 
edges going out of $\mathcal{T}$,
and try to classify such surgeries. 
We focus on Seifert surgeries on torus knots which have 
hyperbolic seiferters or hyperbolic annular pairs. 
By Thurston's hyperbolic Dehn surgery theorem \cite{T1, T2, BP, PetPorti, BoileauPorti}, 
if $(T_{p, q}, m)$ has a hyperbolic seiferter
(resp.\ a hyperbolic annular pair),
then all but finitely many vertices of the 1--complex generated by
successive twists along the seiferter (resp.\ the annular pair)
are Seifert surgeries on hyperbolic knots.
Hence, we call $(T_{p,q}, m)$ with a hyperbolic seiferter or
a hyperbolic annular pair a \textit{spreader}. 
Previously known examples of spreaders 
\cite{DMM2, DMM1, DEMM, DMMtrefoil} have specific patterns
and lead us to the following conjecture.  

\begin{conjecture}
\label{spreaderQ}
If $(T_{p, q}, m)$ is a spreader, 
then $q = 1,\ 2$, or $m= pq,\  pq \pm 1$. 
\end{conjecture}

In Section~\ref{m-move}, we review the definition of $m$--moves
introduced in \cite{DMM1}, which are in fact Kirby calculus
handle--slides over $m$--framed knots.
A trivial knot obtained from a seiferter for $(K, m)$
by a sequence of $m$--moves is also a seiferter for $(K, m)$
if $K$ is nontrivial.
In Sections~\ref{Om} and \ref{Tp2m}, 
we exploit $m$--moves to find seiferters for $(T_{p, q}, m)$
where $q =1, 2$.
Theorems~\ref{seiferter unknot} and \ref{Tp2}
imply the following.

\begin{theorem}
\label{result1}
For each integer $m$, 
$(T_{p, 1}, m) = (O, m)$ and $(T_{p, 2}, m)$ are spreaders.  
In particular, 
$(O, m)$ has infinitely many hyperbolic annular pairs 
as well as infinitely many hyperbolic seiferters. 
\end{theorem}

Regarding seiferters for $(T_{p, q}, m)$ where $q \ge 3$,
we consider two cases according as
$T_{p, q}(m)$ has a unique Seifert fibration up to isotopy or not:
the case when $|m -pq| \ge 2$ and the case when $m = pq, pq \pm 1$.
In the latter case, we prove the theorem below,
which follows from
Propositions~\ref{non degenerate-equivalent},
\ref{non lens-equivalent} and \ref{non lens-equivalent2}.

\begin{theorem}
\label{result2}
Each of $(T_{p, q}, pq)$ $(q \ge 2,  (p, q) \ne (\pm 3, 2))$ and 
$(T_{2n \pm 1, n}, n(2n \pm 1) -1)$ $(n \ge 2)$ has
a hyperbolic seiferter which
cannot be obtained from basic seiferters or a regular fiber of 
the exterior of the torus knot by any sequence of $m$--moves.  
\end{theorem}

Conjecture~\ref{spreaderQ} above
implies that if $q \ge 3$ and $m \ne pq, \ pq \pm 1$,
$(T_{p, q}, m)$ has no hyperbolic seiferters. 
Theorem~\ref{result3} below shows the difficulty of
obtaining such a hyperbolic seiferter.

\begin{theorem}
\label{result3}
Suppose that $q \ge 3$ and $m \ne pq,\  pq \pm 1$
$($i.e.\ $T_{p, q}(m)$ is not a connected sum of
lens spaces, a lens space, or a prism manifold$)$.
Then every seiferter for $(T_{p, q}, m)$ 
is obtained from a basic seiferter or a regular fiber of $S^3 - N(T_{p, q})$ 
by a sequence of $m$--moves $($Proposition~\ref{seiferter for T}$)$.  
However, to obtain a hyperbolic seiferter for $(T_{p, q}, m)$
in such a manner
we need to apply $m$--moves at least twice
$($Corollary~\ref{cor:seiferter for T}(2)$)$.
\end{theorem}

We close the introduction with the following question.

\begin{question}
\label{lens_spreader}
Does every lens surgery $(T_{p, q}, pq \pm 1)$ have a hyperbolic seiferter? 
\end{question}

\noindent
\textbf{Acknowledgments.} 
We would like to thank Mario Eudave-Mu\~noz, Cameron Gordon, John Luecke, 
Daniel Matignon, Jos\'e Mar\'ia Montesinos, Kunio Murasugi, Kouki Taniyama, Masakazu Teragaito, 
and Akira Yasuhara for useful discussions and encouragements. 
We would also like to thank the referee for careful reading and useful comments. 

The third author has been partially supported by JSPS Grants-in-Aid for Scientific 
Research (C) (No.17540097 and No.21540098), 
The Ministry of Education, Culture, Sports, Science and Technology, Japan and 
Joint Research Grant of Institute of Natural Sciences at Nihon University for 2012.

\section{Seiferters for torus knots and Seifert fibrations of torus knot spaces}
\label{m-move}

\begin{figure}[htbp]
\begin{center}
\includegraphics[width=1.0\linewidth]{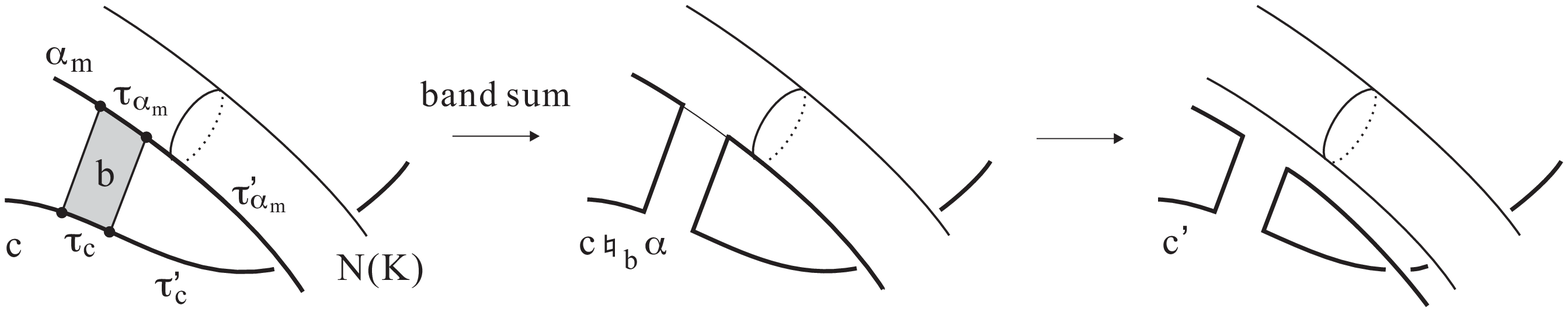}
\caption{$m$--move}
\label{movepush}
\end{center}
\end{figure}

\begin{definition}[\textbf{$m$--move}]
Let $K$ be a knot in $S^3$, and $c$ a knot in $S^3 - N(K)$. 
Take a simple closed curve $\alpha_m$ on $\partial N(K)$ 
representing a slope $m$. 
Let $b$ be a band in $S^3 -\mathrm{int}N(K)$
connecting $\alpha_m$ and $c$,
and let $b \cap \alpha_m = \tau_{\alpha_m}$, $b \cap c = \tau_c$.
We set $\tau'_{\alpha_m} = \alpha_m - \mathrm{int}\tau_{\alpha_m}$
and $\tau'_{c} = c - \mathrm{int}\tau_{c}$. 
Then the band connected sum 
$c\,\natural_b\,\alpha_m
= \tau'_c \cup (\partial b - \mathrm{int}(\tau_c \cup \tau_{\alpha_m}))
\cup \tau'_{\alpha_m}$ is a knot in $S^3 -\mathrm{int}N(K)$.
Pushing $c\,\natural_b\,\alpha_m$ away 
from $\partial N(K)$, 
we obtain a knot $c'$ in $S^3 - N(K)$; 
see Figure~\ref{movepush}.
We say that $c'$ is obtained from $c$ by an \textit{$m$--move}
using the band $b$. 
\end{definition}

\begin{figure}[htbp]
\begin{center}
\includegraphics[width=1.0\linewidth]{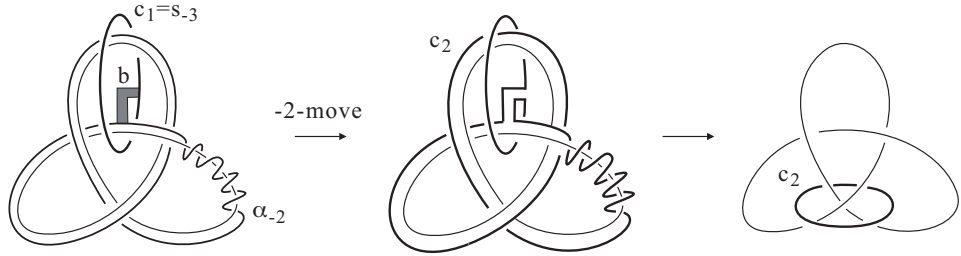}
\caption{$m$--move; $m= -2$, and
$c_2$ is a seiferter for $(T_{-2, 3}, -2)$.}
\label{minus2move}
\end{center}
\end{figure}

Let $K$ be a knot in $S^3$, and $c_1, c_2$ knots in $S^3 - N(K)$.
Assume that $c_2$ is obtained from $c_1$
after a finite sequence of $m$--moves and
isotopies in $S^3 - \mathrm{int}N(K)$.
We then say that
$c_2$ is \textit{$m$--equivalent} to $c_1$. 
Note that $c_2$ is isotopic to $c_1$ in the surgered manifold $K(m)$
(\cite[Proposition~2.19(1)]{DMM1}).
Hence,
if $(K, m)$ is a Seifert surgery,
$c_1$ is a seiferter for $(K, m)$, and
$c_2$ is a trivial knot, then $c_2$ is
a possibly irrelevant seiferter for $(K, m)$. 
Proposition~2.19(3) in \cite{DMM1} shows that $c_2$ is not
irrelevant if $K$ is a nontrivial knot.
Figure~\ref{minus2move} illustrates how an $m$--move works, 
where $K = T_{-3, 2}$, $m = -2$, $c_1 = s_{-3}$.
It follows that $c_2$ is a seiferter for $(T_{-3, 2}, -2)$.
See Section~\ref{Om} for $m$--moves of
annular pairs of seiferters. 

Most seiferters for $(T_{p, q}, m)$ are $m$--equivalent to
basic seiferters or regular fibers of 
Seifert fibrations of $S^3 - \mathrm{int}N(T_{p, q})$. 
Precise statements are as follows.

\begin{proposition}
\label{seiferter for T}
Let $T_{p, q}$ be a nontrivial torus knot, and 
$c$ a seiferter for $(T_{p, q}, m)$, where $m \ne pq$. 
\begin{enumerate}
\item
Suppose that $c$ is an exceptional fiber 
in some Seifert fibration of $T_{p, q}(m)$. 
If $T_{p, q}(m)$ is a lens space, 
we assume that the base surface is $S^2$.  
Then $c$ is $m$--equivalent to a basic seiferter 
$s_p$, $s_q$ or $c_{\mu}$. 
\item
Suppose that $c$ is a regular fiber 
in some Seifert fibration of $T_{p, q}(m)$. 
If $T_{p, q}(m)$ is neither a lens space nor a prism manifold, 
then $c$ is $m$--equivalent to a regular fiber 
in $S^3 - N(T_{p, q})$. 
\end{enumerate}
\end{proposition}

\noindent
\textit{Proof of Proposition~\ref{seiferter for T}.}
We denote by $\mathcal{F}$
a Seifert fibration on $T_{p, q}(m)$ in which $c$ is
an exceptional fiber or a regular fiber. 

\textit{Case $1$. $T_{p, q}(m)$ is not a lens space.} 

By \cite[Proposition~2.8]{DMM1},
if $T_{p,q}(m)$ admits a degenerate Seifert fibration,
then it is either a lens space or
a connected sum of two lens spaces. 
It follows that $\mathcal{F}$ is a non-degenerate Seifert fibration.  
Let $\mathcal{F}_0$ be a natural extension of
the Seifert fibration of $S^3 -\mathrm{int}N( T_{p, q} )$ over $T_{p, q}(m)$. 
The base space of $\mathcal{F}_0$ is the 2--sphere,
and its exceptional fibers are $s_p, s_q$,
and a core of the filled solid torus,
whose indices are $|p|, q$, and $|pq - m|$, respectively.
We note that $c_{\mu}$ is isotopic in $T_{p, q}(m)$
to the third exceptional fiber of $\mathcal{F}_0$. 
Let $t$ be a regular fiber 
of $\mathcal{F}_0|( S^3 - N(T_{p, q}) )$.

\textit{Subcase $1$.
$T_{p, q}(m)$ is not a prism manifold.}
 
It then follows from \cite[Corollary~3.12]{JWW} that 
two Seifert fibrations $\mathcal{F}$ and $\mathcal{F}_0$ are isotopic. 
Hence,
if $c$ is an exceptional fiber in $T_{p, q}(m)$,
then $c$ is isotopic to
one of $s_p, s_q$ and $c_{\mu}$ in $T_{p, q}(m)$,
and thus $m$--equivalent to a basic seiferter 
$s_p, s_q$ or $c_{\mu}$.
Similarly, if $c$ is a regular fiber in $T_{p, q}(m)$,
then $c$ is isotopic to $t$ in $T_{p, q}(m)$
and thus $m$--equivalent to the regular fiber $t$ (\cite[Proposition~2.19(1)]{DMM1}).

\textit{Subcase $2$. 
$T_{p, q}(m)$ is a prism manifold and $c$ is an exceptional fiber.}

A Seifert fibration of a prism manifold is
either over $S^2$ with three exceptional fibers of indices $2, 2, x$ or
over $\mathbb{R}P^2$
with at most one exceptional fiber (\cite[VI.16(b)]{J}).
Hence, $\mathcal{F}_0$ is a Seifert fibration over
the base orbifold $S^2(2, 2, x)$
for some odd integer $x (\ge 3)$. 
Now let us show that $T_{p, q}(m)$ has a Seifert fibration over $S^2$ 
with $c$ an exceptional fiber
even if the base space of $\mathcal{F}$ is not $S^2$.
Assume that $\mathcal{F}$ is a Seifert fibration over $\mathbb{R}P^2$; 
then $\mathcal{F}|( T_{p, q}(m) -\mathrm{int}N(c) )$ is
a Seifert fibration over the M\"obius band
with no exceptional fiber. 
Hence $T_{p, q}(m)- \mathrm{int}N(c)$ admits a Seifert fibration
over the disk with two exceptional fibers of indices $2, 2$.
Extending this fibration over $T_{p,q}(m)$,
we obtain a Seifert fibration over $S^2$
with $c$ an exceptional fiber, as claimed.
For simplicity, denote the new Seifert fibration
by the same symbol $\mathcal{F}$. 
Then, $\mathcal{F}$ is a Seifert fibration over
the base orbifold $S^2(2, 2, x')$ for some odd integer $x'(\ge 3)$. 
Since a regular fiber of $\mathcal{F}$ (resp.\ $\mathcal{F}_0$)
generates the center of $\pi_1( T_{p, q}(m) )$,
the quotient of $\pi_1( T_{p, q}(m) )$
by its center is the dihedral group of order $2x'$ (resp.\ $2x$).
It follows that $x = x'$.

\begin{claim}
\label{prism}
There exists an orientation preserving homeomorphism $f$ 
of $T_{p, q}(m)$ which carries fibers of $\mathcal{F}$ to 
fibers of $\mathcal{F}_0$. 
\end{claim}

\noindent
\textit{Proof of Claim~\ref{prism}.}
We denote the normalized Seifert invariant
of $\mathcal{F}$ by
$(b, \frac{1}{2}, \frac{1}{2}, \frac{y}{x})$ $(b \in \mathbb{Z},\ 0< y < x)$,
and that of $\mathcal{F}_0$ by
$(b', \frac{1}{2}, \frac{1}{2}, \frac{y'}{x})$ $(b' \in \mathbb{Z},\ 0< y' < x)$). 
Note that the order of $H_1(T_{p, q}(m))$ is given by 
$4|(b+1)x + y| = 4|(b'+1)x + y'|$. 
Hence we have $b = b',\ y = y'$ or
$b + b' = -3,\ x = y + y'$. 
In the former case, 
we have an orientation preserving homeomorphism 
of $T_{p, q}(m)$ which carries fibers of $\mathcal{F}$ to 
those of $\mathcal{F}_0$ as desired; 
see \cite{Orlik}, \cite{NR}, and \cite{Hat2}. 
We show that the latter does not occur. 
If we have the latter case, 
then 
$(b', \frac{1}{2}, \frac{1}{2}, \frac{y'}{x}) 
= (-b-3, \frac{1}{2}, \frac{1}{2}, \frac{x-y}{x})$. 
On the other hand, 
$-T_{p, q}(m)$ ($T_{p, q}(m)$ with orientation reversed) 
has a Seifert invariant 
$(-b, -\frac{1}{2}, -\frac{1}{2}, -\frac{y}{x})$, 
which is normalized to 
$(-b-3, \frac{1}{2}, \frac{1}{2}, \frac{x-y}{x})$. 
Thus we have an orientation preserving homeomorphism 
from $-T_{p, q}(m)$ to $T_{p, q}(m)$ 
(\cite{Orlik}, \cite{NR}, and \cite{Hat2}), 
i.e. $T_{p, q}(m)$ admits an orientation reversing homeomorphism. 
This contradicts the fact that a prism manifold has no orientation reversing 
homeomorphism (\cite{Asano}, \cite[8.4]{NR}, \cite{Rubinstein}). 
\QED{Claim~\ref{prism}}

Then, \cite[Lemma~3.5]{JWW} implies that 
$f$ is isotopic to 
a homeomorphism preserving $\mathcal{F}$. 
This implies that $\mathcal{F}_0$ is isotopic to $\mathcal{F}$. 
Hence just as in Subcase~1, 
the exceptional fiber $c$ is $m$--equivalent to one of $s_p, s_q$ and $c_{\mu}$. 
\par

\textit{Case $2$. 
$T_{p, q}(m)$ is a lens space, 
and $c$ is an exceptional fiber.}
  
Then $T_{p, q}(m)$ has a natural Seifert fibration over $S^2$ 
in which $s_p$ and $s_q$ are exceptional fibers of indices 
$|p|, q$. 
Note also that $s_p$ and $s_q$ give a genus one Heegaard splitting
$T_{p, q}(m) =V \cup W$ of 
the lens space $T_{p, q}(m)$; 
$s_p$ and $s_q$ are cores of the solid tori $V$ and $W$.
We recall that the base space of the Seifert fibration $\mathcal{F}$
is $S^2$ from the assumption of Proposition~\ref{seiferter for T}(1).
Then, $\mathcal{F}$ also gives a genus one Heegaard splitting
$T_{p, q}(m) = V' \cup W'$ such that
the exceptional fiber $c$ in $\mathcal{F}$
is a core of $V'$. 
Since a genus one Heegaard splitting
is unique up to isotopy by \cite{Bonahon, HR}, 
$c$ is isotopic to $s_p$ or $s_q$ in $T_{p, q}(m)$.
Proposition~2.19(1) in \cite{DMM1} thus  shows that $c$ is $m$--equivalent to
a basic seiferter $s_p$ or $s_q$ as desired.  
\QED{Proposition~\ref{seiferter for T}}

\begin{remark}
\label{degenerate, lens, prism}
Assumptions in Proposition~\ref{seiferter for T} are necessary.
\begin{enumerate}
\item
As we will see in Proposition~\ref{non degenerate-equivalent}, 
each $(T_{p, q}, pq)$ where $(p, q) \ne (\pm 3, 2)$ has a seiferter which is not $pq$--equivalent to 
any basic seiferter nor a regular fiber of $S^3 - N(T_{p, q})$. 

\item
If $T_{p, q}(m)$ is a prism manifold $($i.e. $q=2$ and $m=2p\pm 2)$, 
then there exists a seiferter $c$ for $(T_{p,q}, m)$
which is a regular fiber in a Seifert fibration over the projective plane (\cite[Corollary 3.15(6)]{DMM1}). 
Then $c$ is not $m$--equivalent to a regular fiber
of $S^3 - N(T_{p, q})$. 

\item
Propositions~\ref{non lens-equivalent} and 
\ref{non lens-equivalent2} show that
for some lens surgeries $(T_{p, q}, m)$ $(m =pq\pm 1)$,
there exist seiferters which are not $m$--equivalent to 
any basic seiferters nor regular fibers of $S^3 - N(T_{p, q})$.
\end{enumerate}
\end{remark}

\section{Annular pairs of seiferters for $(O, m)$}
\label{Om}

Let $\{c_1, c_2\}$ be an annular pair of seiferters.
When we mention the linking number $\mathrm{lk}(c_1, c_2)$,
$c_1$ and $c_2$ are oriented so as to be homologous in
an annulus cobounded by $c_1, c_2$.
If $c_1\cup c_2$ is not a Hopf link,
then this convention determines the linking number
without specifying the annulus.
A Hopf link cobounds two non-isotopic annuli
according as $\mathrm{lk}(c_1, c_2) =1$ or $-1$. 
For details see Lemma~2.30 and Remark~2.31 in \cite{DMM1}.

In \cite{DMM1} an annular pair $\{c_1, c_2\}$ is defined
to be an ordered pair of $c_1$ and $c_2$ to specify
the direction of twist along the annulus
cobounded by $c_1 \cup c_2$. 
However, since we do not perform annulus twists in this paper,
annular pairs are presented as unordered pairs.

Let $K$ be a knot in $S^3$, and $c_1 \cup c_2$ a link in $S^3-N(K)$.
Let $c'_1$ be a knot obtained from $c_1$ by an $m$--move
using a band disjoint from $c_2$
and connects $c_1$ and a simple closed curve
on $\partial N(K)$ with slope $m$.
We then say that
$c'_1 \cup c_2$ is obtained from $c_1 \cup c_2$
by an $m$--move.
The link $c'_1 \cup c_2$ is isotopic to $c_1\cup c_2$ in
the surgered manifold $K(m)$ as ordered links
(\cite[Lemma~2.25(1)]{DMM1}).
If $\{c_1, c_2\}$ is a pair of seiferters for a Seifert surgery
$(K, m)$ and $c'_1$ is trivial in $S^3$, then
$\{c'_1, c_2\}$ is also a pair of seiferters for $(K, m)$ 
(\cite[Lemma~2.25(2)]{DMM1}).
The theorem below complements Theorem~6.21 in \cite{DMM1}.  

\begin{theorem}
\label{seiferter unknot}
\begin{enumerate}
\item 
For each integer $m$, 
there are infinitely many hyperbolic Hopf pairs of seiferters
for $(O, m)$.
\item 
For any integers $m \ne 0$ and $p \ge 2$ except 
$(m, p) = (\pm 1, 2)$, 
there is a hyperbolic annular pair of seiferters 
$\{ c_1, c_2 \}$ for $(O, m)$ with $\mathrm{lk}(c_1, c_2) = p$. 
\end{enumerate}
\end{theorem}

\noindent
\textit{Proof of Theorem~\ref{seiferter unknot}.}
(1) Assertion~(1) follows from Lemma~\ref{ccpHopf} below. 
\par
(2) Assume that $m \ne 0$, $p \ge 2$, and $(m, p) = (\pm 1, 2)$.
Then, if $m \ne 1$, $\{ c, c_{p+1, m} \}$ in Proposition~\ref{hyp ann pair1}
with $q$ replaced by $p+1$ is a hyperbolic annular pair for $(O, m)$ with 
$\mathrm{lk}(c, c_{p+1, m}) =p$, as desired in assertion~(2).
If $m \ne -1$, $\{ c, c'_{p-1, m} \}$ in Proposition~\ref{hyp ann pair2}
has the desired property.
\QED{Theorem~\ref{seiferter unknot}}

\begin{figure}[hbt]
\begin{center}
\includegraphics[width=0.8\linewidth]{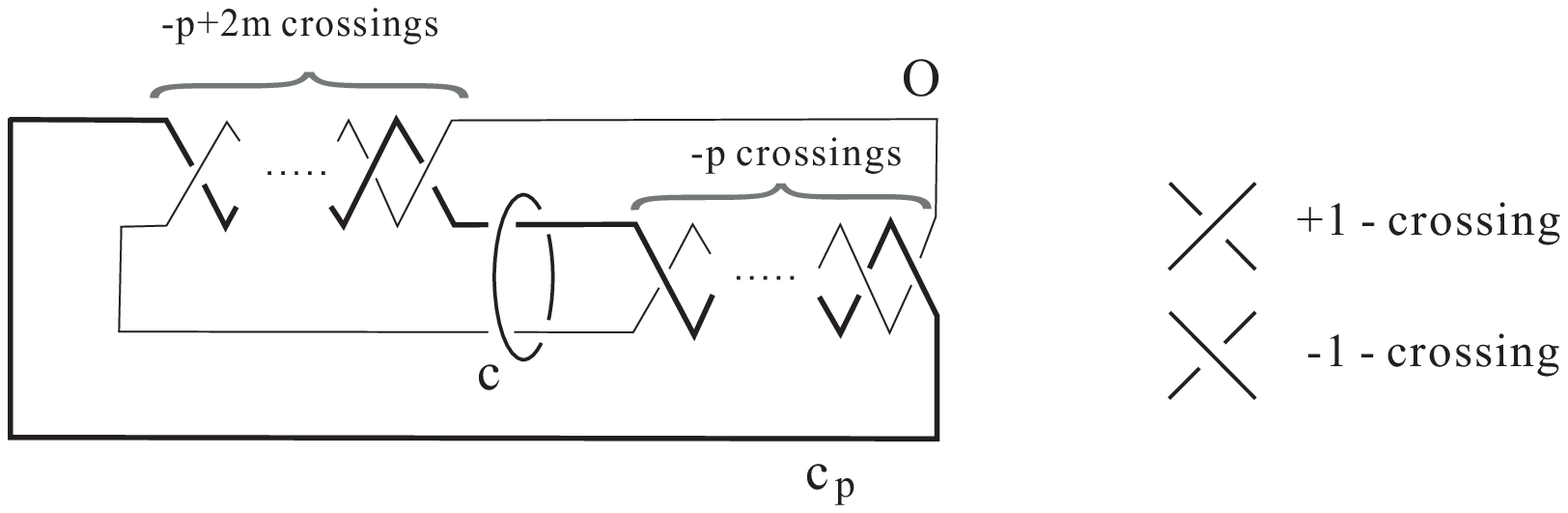}
\caption{}
\label{Omseiferters}
\end{center}
\end{figure}

\begin{lemma}
\label{ccpHopf}
Let $O \cup c \cup c_p$ be the link in Figure~\ref{Omseiferters},
where $p$ is an odd integer with $|p| \ge 3$. 
Then, $\{c, c_p\}$ is a hyperbolic Hopf pair of seiferters
for $(O, m)$ if $p \ne 2m \pm 1$. 
For each $m$,
$\{c, c_p\}$ $(p \ge m, p\ne 2m \pm 1)$ are mutually distinct,
hyperbolic Hopf pairs.
\end{lemma}

\noindent
\textit{Proof of Lemma~\ref{ccpHopf}.}
Consider the link consisting of
a torus knot $T_{p, 2}$ ($|p| \ge 3$) and
its basic seiferters $s_2$, $s_p$.
Regard $s_2$ as the trivial knot $O$, and set $c =s_p$.
See the first figure of Figure~\ref{Omseiferters2}.
There is a Seifert fibration of $S^3 -\mathrm{int}N(O)$
in which
$T_{p, 2}$ is a regular fiber and $c$ is the exceptional fiber
of index $|p|$.
Let $c_p$ be the trivial knot obtained from $T_{p, 2}$ in $S^3 -N(O)$ 
by the $m$--move in Figure~\ref{Omseiferters2}.
Since $c \cup T_{p, 2}$ is isotopic in $O(m)$ to $c \cup c_p$, 
$c \cup c_p$ is also the union of fibers in a Seifert fibration
of $O(m)$. 
It follows that $\{c, c_p\}$
in the second figure of Figure~\ref{Omseiferters2}
is a pair of seiferters for $(O, m)$.
After isotopy, the link $O \cup c \cup c_p$
in the last figure of Figure~\ref{Omseiferters2}
is the same link as $O \cup c \cup c_p$ in Figure~\ref{Omseiferters}.
Hence, $\{c, c_p\}$ in Figure~\ref{Omseiferters} is
a Hopf pair of seiferters for $(O, m)$.

\begin{figure}[hbt]
\begin{center}
\includegraphics[width=1.0\linewidth]{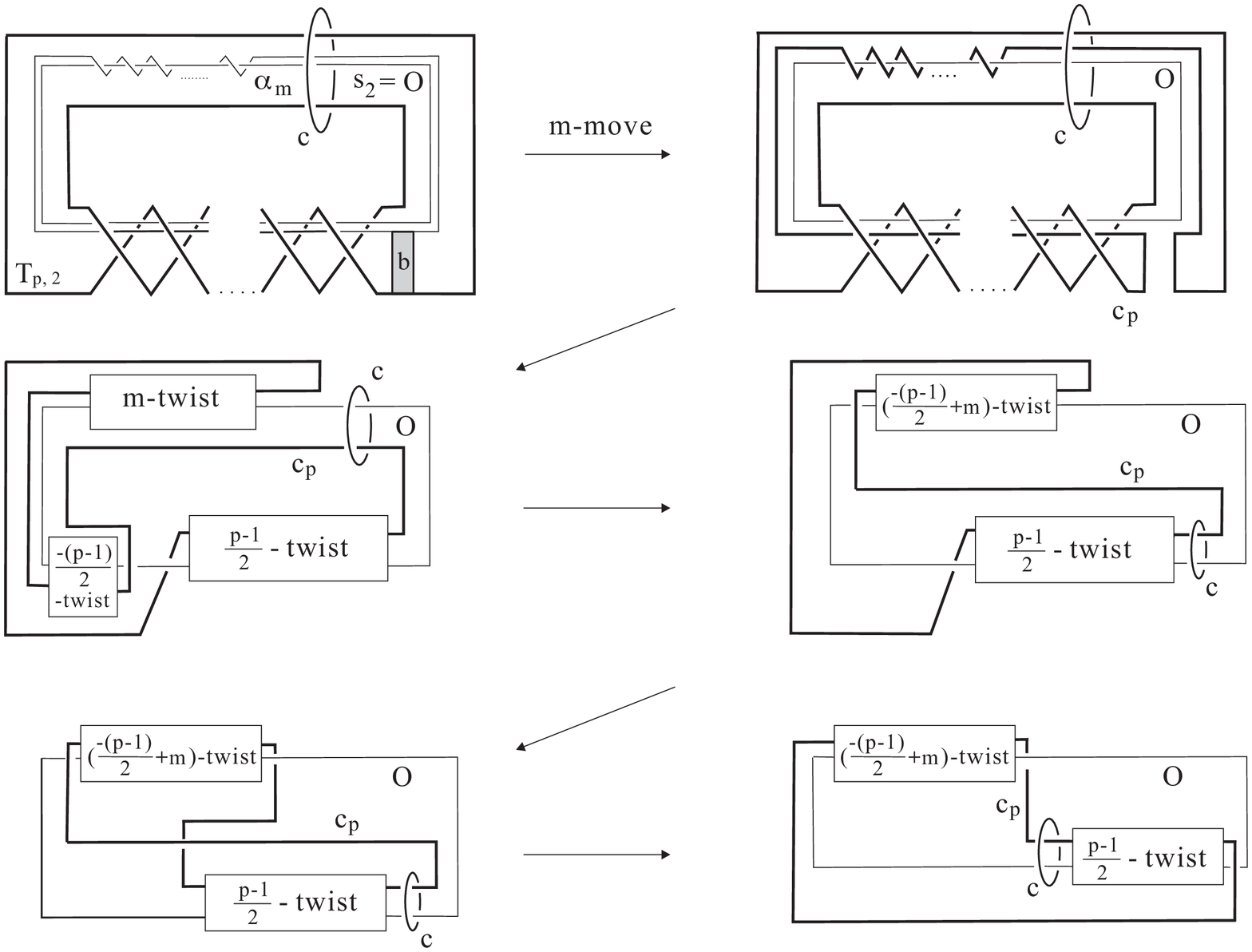}
\caption{}
\label{Omseiferters2}
\end{center}
\end{figure}

Suppose $p \ne 2m \pm 1$.
Let us show that $\{c, c_p\}$ is a hyperbolic Hopf pair,
i.e.\ $O \cup c \cup c_p$ is a hyperbolic link. 
Assume for a contradiction that
$X = S^3 -\mathrm{int}N(O \cup c \cup c_p)$ is Seifert fibered.
Then, the exterior of $O \cup c_p$,
which is obtained from $X$ by Dehn filling along $\partial N(c)$,
is a non-degenerate Seifert fiber space or a reducible manifold.
On the other hand, since
$O \cup c_p$ $(p \ne 2m \pm 1)$ is a $2$--bridge link
and not a torus link,
it is a hyperbolic link.
(For details refer to the proof of Theorem~6.21 in \cite{DMM1}.)
This is a contradiction,
so that $X$ is not Seifert fibered.
Figure~\ref{Omseiferters3} shows that
$X$ is homeomorphic to the exterior of
the Montesinos link $L = M(\frac{1}{2m-p-1}, \frac{1}{2}, \frac{1}{2})$. 
The proof of \cite[Corollary~5]{Oertel} shows
that $X$ is hyperbolic if $X$ is not Seifert fibered, and $L$ is not
equivalent to the Montesinos links
$M(\frac{1}{2}, \frac{1}{2}, \frac{-1}{2}, \frac{-1}{2})$,
$M(\frac{2}{3}, \frac{-1}{3}, \frac{-1}{3})$,
$M(\frac{1}{2}, \frac{-1}{4}, \frac{-1}{4})$, 
$M(\frac{1}{2}, \frac{-1}{3}, \frac{-1}{6})$,
or the mirror images of these links.
The 2--fold branched cover of $S^3$ along $L$ is a prism manifold, 
which has a finite fundamental group.
However, the 2--fold branched covers along 
the four Montesinos links above have infinite fundamental groups. 
Therefore, $X$ is hyperbolic.

We note that
$\{ |\mathrm{lk}(c, O)|,\, |\mathrm{lk}(c_p, O)| \}
= \{1, |m-p|\}$.
Hence, if $O \cup c \cup c_p$ is isotopic to
$O \cup c \cup c_q$ $(p, q \ge m)$ in $S^3$ with $O$ sent to $O$,
then $p=q$.
It follows that
for each $m$ the pairs of seiferters $\{c, c_p\}$ where $p\ge m$ are
mutually distinct.
\QED{Lemma~\ref{ccpHopf}}

\begin{figure}[hbt]
\begin{center}
\includegraphics[width=1.0\linewidth]{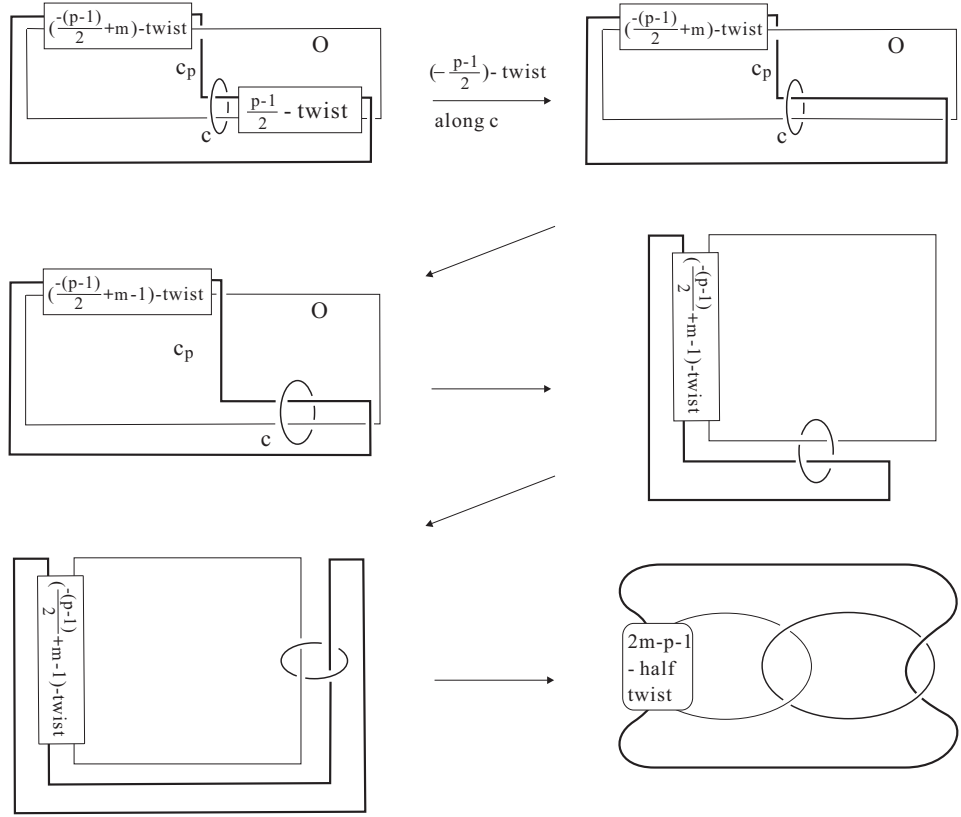}
\caption{Continued from Figure~\ref{Omseiferters2}.}
\label{Omseiferters3}
\end{center}
\end{figure}

\begin{remark}
\label{error:oertel}
Corollary~5 in \cite{Oertel} states that
a Montesinos link is hyperbolic if
it is not a torus link, and not equivalent to the four Montesinos links
listed above or their mirror images.
However, in the proof the author assumes that
links whose exteriors are Seifert fibered are torus links,
which is not true.
We thus obtain the corrected Corollary~5 in \cite{Oertel}
by replacing the word ``torus link" with
``link whose exterior is Seifert fibered".
\end{remark}

The Hopf pair of seiferters $\{c, c_p\}$ satisfies
$|\mathrm{lk}(c, c_p)| = 1$. 
Now for a given integer $p>1$,
let us find an annular pair of seiferters $\{c_1, c_2\}$ 
with $\mathrm{lk}(c_1, c_2) = p$
as claimed in Theorem~\ref{seiferter unknot}(2). 
We will give such examples in
Propositions~\ref{hyp ann pair1} and~\ref{hyp ann pair2}. 
To prove the hyperbolicity of these examples, 
we prepare some general results. 

\begin{proposition}
\label{torus bounds solid torus}
Let $l_1 \cup \cdots \cup l_n$ be an $n$-component link in a solid torus $V$.
Suppose that there is a meridional disk $D$ for $V$ satisfying $(1)$, $(2)$ below.

\begin{enumerate}
\item The winding number of $l_i$ in $V$ equals $| D\cap l_i |$
for any $i$.
\item 
$V - \mathrm{int}N(D \cup (\cup_{i=1}^n l_i) )$ is homeomorphic  to a handlebody.
\end{enumerate}

Then, if $\displaystyle V -\mathrm{int}N( \cup_{i=1}^n l_i )$ contains
an essential torus, it bounds a solid torus in $V$.
\end{proposition}

\noindent
\textit{Proof of Proposition}~\ref{torus bounds solid torus}.
We identify $V$ split along $D$ with $D^2\times I$,
where $I =[0, 1]$, and $D^2\times \{0\}$ and $D^2 \times \{1\}$ are
identified with $D$ in $V$.
Let $a_1, \ldots, a_m$ be the arcs in $D^2 \times I$ obtained by cutting
$\cup_{i=1}^n l_i$ by $D$;
each $a_i$ connects $D^2\times \{0 \}$ and
$D^2 \times \{1 \}$ by condition~(1).

Assume that $V -\mathrm{int}N(\cup_{i=1}^n l_i)$ contains an essential torus $T$.
Isotope $T$ in $V -\mathrm{int}N(\cup_{i=1}^n l_i)$ so as to minimize
the number of components $|D \cap T|$.
Note that condition~(2) implies $D \cap T \ne \emptyset$.
Then $D$ splits $T$ into essential annuli $A_1, A_2, \ldots, A_k$
properly embedded in $D^2 \times I -\cup_{i=1}^m a_i$
such that a component of $\partial A_i$ 
and a component of $\partial A_{i +1}$ are identified in $V$,
where $1 \le i \le k$, and if $i =k$ we regard $i+1 =k+1$ as $1$.

\begin{claim}
\label{annulus connects}
Each annulus $A_i$ connects $D^2\times \{ 0 \}$
and $D^2 \times \{1 \}$.
\end{claim}

\noindent
\textit{Proof of Claim}~\ref{annulus connects}.
Assume for a contradiction that
some $A_{i_0}$ satisfies
$\partial A_{i_0} \subset D^2\times \{\alpha\}$,
where $\alpha= 0$ or $1$.
Let $B_1, B_2$ be the disks in $D^2 \times \{\alpha\}$ bounded by
the components of $\partial A_{i_0}$.
If $B_1 \cap B_2 = \emptyset$,
then $B_1$ intersects some arc $a_{j_0}$ and
$B_1 \cup A_{i_0} \cup B_2$ bounds a 3-ball in $D^2 \times I$.
This implies $\partial a_{j_0} \subset D^2 \times \{\alpha \}$,
a contradiction to condition~(1)
in Proposition~\ref{torus bounds solid torus}.
It follows that $B_1 \subset \mathrm{int}B_2$
or $B_2 \subset \mathrm{int}B_1$. 
Without loss of generality, we assume that the former holds.
Let $M$ be the $3$--submanifold in $D^2\times I$
bounded by the torus $A_{i_0}\cup ( B_2 -\mathrm{int}B_1)$.
Condition~(1) then implies that $M \cap a_i = \emptyset$
for any $i$.

\textit{Case $1$. $M$ is boundary irreducible.}
\par
If $\partial M$ is incompressible in
$X = V - \mathrm{int}N(\cup_{i=1}^n l_i) -\mathrm{int}M$,
then after pushing $\partial M$ in $V -\mathrm{int}N(\cup_{i=1}^n l_i)$
off $D$,
$\partial M$ is an essential torus in
$D^2 \times I - \mathrm{int}N( \cup_{i=1}^m a_i )$.
This contradicts condition~(2) in Proposition~\ref{torus bounds solid torus}.
Hence, an essential simple closed curve $c$ in $\partial M$ bounds
a disk in $X$.
On the other hand, 
$\partial B_1$ is also an essential simple closed curve in 
$\partial M$ bounding the disk $B_1$ in $V -\mathrm{int}M$.
Since the rank of $\mathrm{Ker}( H_1( \partial M ) \to H_1( V -\mathrm{int}M ) )$ 
is less than or equal to one by the Poincar\'e duality,
we see that $[c] = [\partial B_1]$ in $H_1( \partial M )$
and thus $\partial B_1$ bounds a disk in $X$.
This contradicts the fact that $A_{i_0}$ is essential in
$D^2 \times I -\cup_{i=1}^m a_i$.

\textit{Case $2$. $M$ is boundary reducible.}
\par
It follows that $M$ is a solid torus.
Since $\partial B_1 (\subset \partial M)$ bounds the disk $B_1$
in $S^3 -\mathrm{int}M$,
a meridian of $M$ and $\partial B_1$ intersect in one point.
This implies that the annulus $A_{i_0}$ is parallel to
$B_1 -\mathrm{int}B_2$ in $M$,
and contradicts the fact that $A_{i_0}$ is essential in
$D^2 \times I -\cup_{i=1}^m a_i$.
\QED{Claim~\ref{annulus connects}}

By Claim~\ref{annulus connects} the union of $A_i$ and the two disks
in $D^2 \times \{0, 1\}$
bounded by $\partial A_i$ bounds a 3--ball $V_i$
in $D^2 \times I$.
Note that for any distinct $i, j$ we have $V_i \cap V_j = \emptyset$,
$V_i \subset V_j -A_j$, or $V_j \subset V_i -A_i$.
If $V_1, V_2, \ldots, V_k$ are mutually disjoint,
then $\cup_{i=1}^k V_i$ forms a solid torus in $V$ bounded by $T$
as claimed in Proposition~\ref{torus bounds solid torus}.
So assume that $V_i \subset V_j -A_j$ for some $i, j$.
Then by Claim~\ref{annulus connects},
$V_{i +\varepsilon} \subset
V_{j +\varepsilon} -A_{j + \varepsilon}$,
where $\varepsilon= \pm 1$ and
we regard $k+1$, $0$ as $1$, $k$, respectively.
Repeating this argument, we see that for any $V_i$ there exists $V_j$
such that $V_i \subset V_j -A_j$.
This does not occur for a finite number of 3-balls
$V_1, V_2, \ldots, V_k$.
This completes the proof of Proposition~\ref{torus bounds solid torus}.
\QED{Proposition~\ref{torus bounds solid torus}}

The following proposition will be useful. 

\begin{proposition}
\label{hyp ann pair}
Let $l_1 \cup l_2$ be a $2$-component link in a solid torus $V$
such that $l_1$ is a $(p, q)$ cable of $V$ where $q \ge 2$,
$l_2$ is a core of $V$, and $l_1 \cup l_2$ satisfies
conditions~$(1)$, $(2)$ in Proposition~$\ref{torus bounds solid torus}$.
Then, $l_1 \cup l_2$ is a hyperbolic link in $V$
if we cannot isotope $l_2$ in $V -\mathrm{int}N(l_1)$ so as to be
disjoint from a cabling annulus for $N(l_1) \subset V$.
\end{proposition}

\noindent
\textit{Proof of Proposition~\ref{hyp ann pair}.}
First we remark that since $l_i$ $(i=1, 2)$ wraps $V$ geometrically at least once, 
$V -\mathrm{int}N( l_1 \cup l_2 )$ is irreducible. 

Assume for a contradiction that $V -\mathrm{int}N( l_1 \cup l_2 )$
contains an essential torus $T$.

\begin{claim}
\label{parallel torus}
$T$ is parallel to $\partial V$ and separates $l_1$ and $l_2$,
and $l_2$ lies between $T$ and $\partial V$.
\end{claim}

\noindent
\textit{Proof of Claim~\ref{parallel torus}.}
Since $T$ is not essential in $V -\mathrm{int}N(l_1)$,
there are three cases: (1)~$T$ is compressible
in $V -\mathrm{int}N(l_1)$,
(2)~$T$ is parallel to $\partial N(l_1)$ in $V- \mathrm{int}N(l_1)$,
and
(3)~$T$ is parallel to $\partial V$ in $V- \mathrm{int}N(l_1)$.
Case~(3) implies Claim~\ref{parallel torus}.
So we derive a contradiction in cases~(1), (2).

Let $V'$ be the solid torus in $V$ bounded by $T$
(Proposition~\ref{torus bounds solid torus});
$V'$ contains at least one of $l_1$ and $l_2$.
Since each $l_i$ is not contained in a $3$--ball in $V$,
$V'$ is not contained in a $3$--ball in $V$, either.
It follows that $T = \partial V'$ is incompressible
in $V -\mathrm{int}V'$. 
Now assume case~(1) occurs.
Then $l_2 \subset V'$,
and $T$ separates $l_1$ and $l_2$.
Since $V -\mathrm{int}N(l_2) \cong T^2 \times I$,
$T$ is parallel to $\partial N(l_2)$
in $V -\mathrm{int}N(l_1 \cup l_2)$.
This contradicts the fact that $T$ is essential.
Assume case~(2) occurs.
Since $T$ is essential in $V-\mathrm{int}N(l_1 \cup l_2)$,
we have $l_2 \subset V'$.
Then the winding number of $l_2$ in $V$
is a multiple of $q( \ge 2)$.
This contradicts the fact that
$l_2$ is a core of $V$.
\QED{Claim~\ref{parallel torus}}

\begin{claim}
\label{off cabling annulus}
For any cabling annulus $A$ for $N(l_1)$ in $V$,
we can isotope $l_2$ in $V -\mathrm{int}N( l_1 )$
so as to be disjoint from $A$.
\end{claim}

\noindent
\textit{Proof of Claim~\ref{off cabling annulus}.}
Let $W$ be the submanifold of $V -\mathrm{int}N( l_1 )$ cobounded by
$\partial V$ and $T$.
Identify $W$ and $\partial V \times I$ so that
$\partial V$ and $T$ correspond to $\partial V \times \{0 \}$
and $\partial V \times \{ 1 \}$ respectively,
and let $\pi: W = \partial V \times I \to I$ be the natural projection.
Since $A' = A \cap W$ is a compact submanifold of $W$,
$\pi( A' )$ is a compact and thus closed subset of $I$.
It follows that $\inf \pi(A') \in \pi( A' )$
and $0 < \inf \pi(A')$,
since $A \subset \mathrm{int}V$.
Now isotope $l_2 \cup T$ in $W$ so that $T$ becomes
$\partial V \times \{ \frac{1}{2}\inf \pi( A' ) \}$.
After this isotopy $l_2$ becomes disjoint from $A$.
\QED{Claim~\ref{off cabling annulus}}

Claim~\ref{off cabling annulus} contradicts the assumption
in Proposition~\ref{hyp ann pair}.
Hence, 
$V -\mathrm{int}N(l_1 \cup l_2)$ contains no essential torus.

\begin{claim}
\label{no annulus}
$X = V -\mathrm{int}N( l_1 \cup l_2 )$ contains no essential annulus.
\end{claim}

\noindent
\textit{Proof of Claim~\ref{no annulus}.}
Assume for a contradiction that $X$ contains an essential annulus.
Since $X$ contains no essential torus,
and is irreducible and boundary irreducible, 
this assumption implies that $X$ is a Seifert fiber space.
Then $X$ contains an essential annulus $A$ connecting
$\partial N( l_1 )$ and $\partial V$;
note that $A$ is also an essential annulus in the cable space
$V -\mathrm{int}N( l_1 )$.
Take a regular neighborhood $N( \partial V \cup A )$ in $X$.
Then the closure of $\partial N( \partial V \cup A ) -\partial X$ is
a cabling annulus for $N( l_1 )$ in $V$.
Since the cabling annulus is disjoint from $l_2$,
this fact contradicts the assumption in Proposition~\ref{hyp ann pair}.
\QED{Claim~\ref{no annulus}}

The proof of Proposition~\ref{hyp ann pair} is thus completed.
\QED{Proposition~\ref{hyp ann pair}}

\begin{figure}[hbt]
\begin{center}
\includegraphics[width=0.9\linewidth]{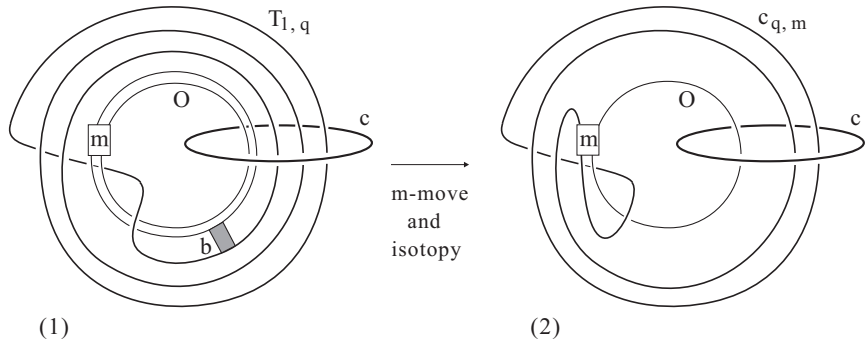}
\caption{Annular pair of seiferters $\{c, c_{q, m}\}$;\ $q= 3$}
\label{Omhighlinking1}
\end{center}
\end{figure}

\begin{proposition}
\label{hyp ann pair1}
Let $c \cup c_{q, m}$ be the link obtained from $c \cup T_{1, q}$ in $S^3 -N(O)$ 
by an $m$-move using the band $b$ in Figure~\ref{Omhighlinking1}(1) and 
an isotopy.  
Assume that $q \ge 3$, $m \ne 0, 1$, and $(m, q) \ne (-1, 3)$.
Then, $\{ c, c_{q, m} \}$ is
a hyperbolic annular pair of seiferters for $(O, m)$ with $\mathrm{lk}(c, c_{q, m}) = q-1$. 
\end{proposition}

\begin{figure}[hbt]
\begin{center}
\includegraphics[width=0.9\linewidth]{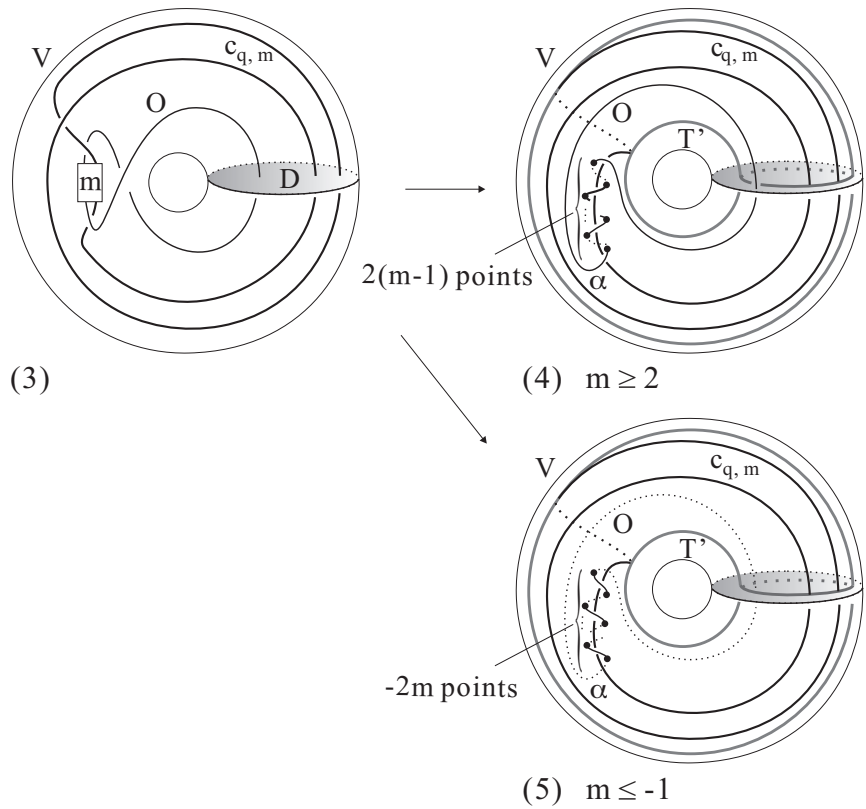}
\caption{In (4), (5),
the intersection points between $O$ and $T'$ are indicated  by ``dots''.}
\label{Omhighlinking2}
\end{center}
\end{figure}

\noindent
\textit{Proof of Proposition $\ref{hyp ann pair1}$}.
Since $S^3 -\mathrm{int}N(O)$ admits a Seifert fibration
in which $c$ and $T_{1, q}$ in Figure~\ref{Omhighlinking1}(1)
are fibers,
$\{c, T_{1, q}\}$ is an annular pair of seiferters for $(O, m)$. 
After the $m$--move in Figure~\ref{Omhighlinking1},
$c_{q, m}$ and $c$ in Figure~\ref{Omhighlinking1}(2) remain fibers in $O(m)$. 
Note that $c_{q, m}$ is the torus knot $T_{1, q -1}$,
a trivial knot in $S^3$,
and $c \cup c_{q, m}$ bounds an annulus.
It follows that
$\{c, c_{q, m}\}$ is an annular pair of seiferters for $(O, m)$. 
Note that $\mathrm{lk}(c, c_{q, m}) = q-1$. 
It remains to show that $O \cup c \cup c_{q, m}$ is a hyperbolic link.

Let $V$ be the solid torus $S^3 -\mathrm{int}N(c)$
containing $O \cup c_{q, m}$.
Then $O$ is a core of $V$ and
$c_{q, m}$ is a $(1, q-1)$ cable of $V$.
The meridional disk $D$ for $V$ described in Figure~\ref{Omhighlinking2}(3)
intersects $O$ in one point and
$c_{q, m}$ in $q-1$ points.
Note also that
$V -\mathrm{int}N(D \cup O \cup c_{q, m})$ is homeomorphic to
a handlebody.
The link $O \cup c_{q, m}$ in $V$ thus satisfies conditions~(1), (2)
with $n =2$ in Proposition~\ref{torus bounds solid torus}.
Then by Proposition~\ref{hyp ann pair},
in order to show that $O \cup c_{q, m}$ is hyperbolic in $V$, 
it is sufficient to show that
$O$ cannot be isotoped in $V -\mathrm{int}N(c_{q, m})$
off a cabling annulus for $N(c_{q, m}) \subset V$.

Now let $T'$ be the torus in $V$ containing $c_{q, m}$
as described in Figure~\ref{Omhighlinking2}(4), (5);
then $A = T' -\mathrm{int}N(c_{q, m})$ is a cabling annulus
for $N(c_{q, m}) \subset V$.
We note that
$T'$ intersects $O$ in $2(m -1)$ points if $m \ge 2$,
and in $-2m$ points if $m \le -1$. 
We denote by $\alpha$ the closure of the component of $O -T'$ intersecting $D$. 
Let $V'$ be the solid torus in $V$ bounded by $T'$.
Concerning the arc components of $O \cap V'$
and $O \cap ( V -\mathrm{int}V' )$,
we can check the following.

\begin{claim}
\label{arcs}
\begin{enumerate}
\item
In $V -\mathrm{int}V'$ $($resp.\ $V'$$)$,
the arc $\alpha$ in Figure~\ref{Omhighlinking2}(4) $($resp.\ (5)$)$
is isotopic with $\partial \alpha$ fixed to an arc in $T'$
intersecting $c_{q, m}$ algebraically twice.
\item
In $V -\mathrm{int}V'$ $($resp.\ $V'$$)$,
each component $\beta$  of $O \cap (V -\mathrm{int}V')$
$($resp.\ $O \cap V'$$)$ other than $\alpha$ is isotopic
with $\partial \beta$ fixed to an arc in $T'$ intersecting
$c_{q, m}$ once.
\end{enumerate}
\end{claim}

Using Claim~\ref{arcs}, we show that
there is no isotopy of $O$ in $V- \mathrm{int}N( c_{q, m} )$
which makes the intersection between $O$ and
the cabling annulus $A$ empty.
Assume for a contradiction
that there is an isotopy $f: S^1\times I \to
V -\mathrm{int}N( c_{q, m} )$
such that $f( S^1\times \{0 \}) =O$
and $f( S^1\times \{1\}) \cap A = \emptyset$.
We may assume that $f$ is transverse to $A$;
then $f^{-1}(A)$ is a $1$-submanifold properly embedded in
$S^1 \times I$.
Since $f( S^1\times \{1\}) \cap A = \emptyset$,
we see $f^{-1}(A) \cap (S^1 \times \{1\}) = \emptyset$,
so that
each arc component of $f^{-1}(A)$
has its end points in $S^1 \times \{ 0\}$. 
If $f^{-1}(A)$ has a circle component bounding a disk in $S^1 \times I$, 
then by the loop theorem and the incompressibility of $A$
in $V -\mathrm{int}N( c_{q, m} )$ $f$ restricted on
the innermost circle is
null-homotopic in $A$. 
Hence we can modify $f$ so that the innermost circle is eliminated. 
Thus by re-choosing $f$ we may assume
$f^{-1}(A)$ does not contain null-homotopic circles
in $S^1 \times I$.
For two arc components $a_1, a_2$ of $f^{-1}(A)$,
we say that $a_1$ is closer to $S^1 \times \{0\}$ than $a_2$
if the disk cobounded by $a_2$ and an arc in $S^1 \times \{0\}$
contains $a_1$.
Let $c_1$ be an arc component of $f^{-1}(A)$
closest to $S^1 \times \{0\}$,
and $c_2$ the arc in $S^1 \times \{0 \}$ such that
$c_1 \cup c_2$ cobounds a disk in $S^1 \times I$.
Note that $f(c_2)$ is the closure of a component of $O - A$,
and $f(c_1)$ is an immersed arc in $A$ with $\partial f(c_1)
= \partial f(c_2)$.

\begin{claim}
\label{q=3}
It holds that $q = 3$, $m \le -1$, and $f(c_2)$ is the arc
$\alpha ( \subset V' )$ in Figure~\ref{Omhighlinking2}$(5)$.
\end{claim}

\noindent
\textit{Proof of Claim~\ref{q=3}.}
Set $X = V'$ if $f(c_2) \subset V'$,
and $X= V -\mathrm{int}V'$ if $f(c_2) \subset V -\mathrm{int}V'$.
Then $f(c_1) (\subset A)$ is homotopic in $X$
to the component $f(c_2)$ of $O \cap X$ with its end points fixed.
Combining this homotopy and the isotopies in Claim~\ref{arcs},
we see that $f(c_1)$ is homotopic in $X$ with its end points fixed
to an arc $\gamma$ in $T'$ intersecting $c_{q, m}$
once (if $f(c_2)$ is an arc $\beta$ in Claim~\ref{arcs}(2))
or algebraically twice (if $f(c_2)$ is the arc $\alpha$ in Claim~\ref{arcs}(1)).
Hence, the closed curve $f(c_1) \cup \gamma$ in $T'$
intersecting $c_{q, m}$ once or algebraically twice
is null-homotopic in $X$.
Since $V -\mathrm{int}V' \cong T^2 \times I$,
$f(c_1) \cup \gamma$, which is not null-homotopic in $T'$,
is not null-homotopic in $V -\mathrm{int}V'$.
It follows that $X = V'$ and thus $f(c_2) \subset V'$.
Since $c_{q, m}$ is the $(1, q-1)$ cable of $V'$ where $q \ge 3$,
a meridian of $V'$ intersects $c_{q, m}$
algebraically $q -1$ times.
It follows that $q =3$ and
$\gamma$ intersects $c_{q, m}$ algebraically twice.
Furthermore, we see that
$f(c_2)$ is the arc $\alpha$ in Figure~\ref{Omhighlinking2}(5) and so $m \le -1$.
\QED{Claim~\ref{q=3}}

\begin{figure}[hbt]
\begin{center}
\includegraphics[width=0.35\linewidth]{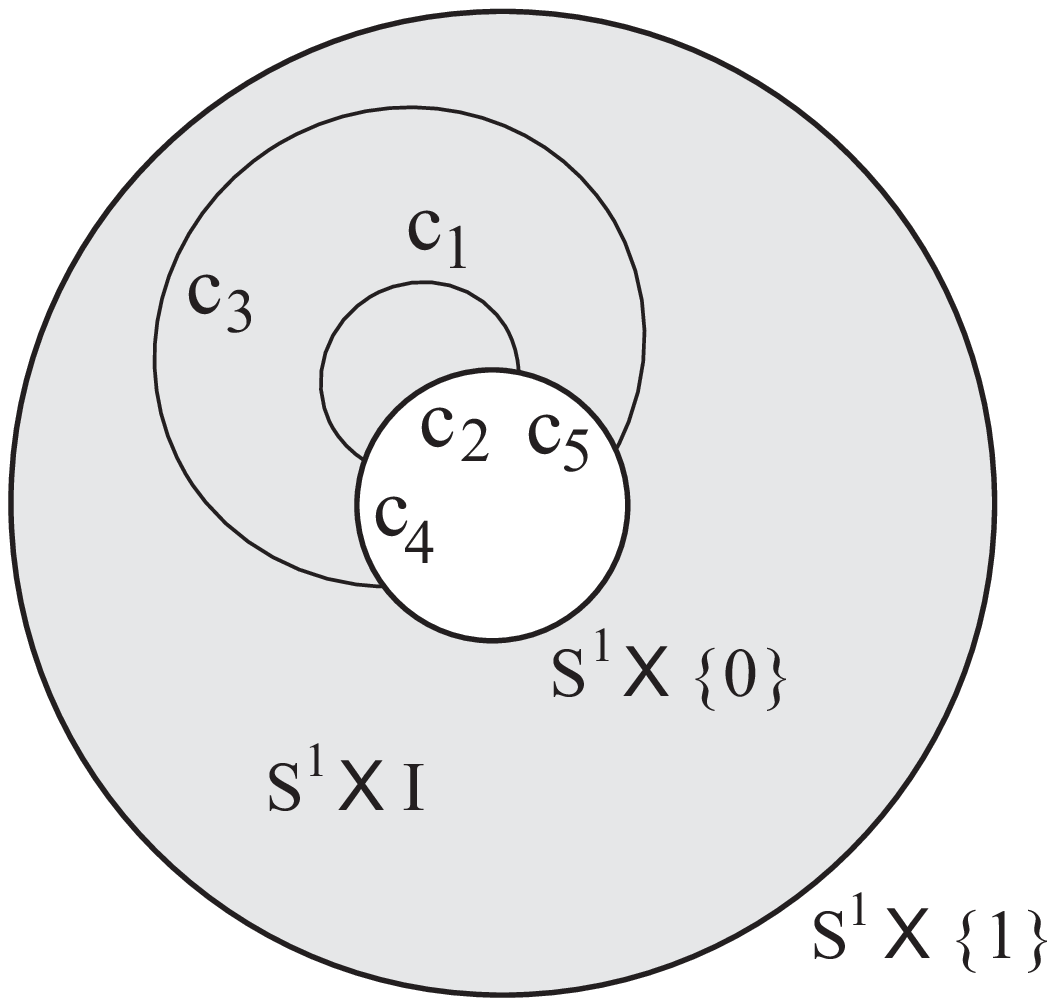}
\caption{}
\label{Omhighlinking3}
\end{center}
\end{figure}

By Claim~\ref{q=3}
$c_1$ is the only arc component of $f^{-1}(A)$ closest to
$S^1 \times \{0\}$.
Hence all arc components of $f^{-1}(A)$ are
parallel to $c_1$ in $S^1 \times I$.
The assumption $(m, q) \ne (-1, 3)$ in 
Proposition~\ref{hyp ann pair1}
together with Claim~\ref{q=3} implies 
$m \le -2$,
so that $A$ intersects $O$ in $-2m ( \ge 4 )$ points and hence 
$f^{-1}(A)$ has
at least two arc components.
Let $c_3$ be the second closest arc component of $f^{-1}(A)$
to $S^1 \times \{ 0\}$,
and $c_4, c_5$ the subarcs of $S^1\times \{ 0\}$
connecting $\partial c_1$ and $\partial c_3$;
see Figure~\ref{Omhighlinking3}.
Note that
$f(c_4)$ and $f(c_5)$ are the components of $O \cap (V -\mathrm{int}V')$
adjacent to $f(c_2) = \alpha$.
Now we give the arcs $c_4$ and $c_5$ the orientations induced
from an orientation of $S^1 \times \{0\}$.  
Then, Figure~\ref{Omhighlinking2}(5) shows that
$f(c_4)$ and $f(c_5)$ are isotopic to arcs in $T'$
whose algebraic intersection numbers with $c_{q, m}$ are both one
under an adequate orientation of $c_{q, m}$.
This implies that
the closed curve $f(c)$ where $c = c_1 \cup c_5 \cup c_3 \cup c_4$
is homotopic in $V -\mathrm{int}V'$ to a closed curve in $T'$
intersecting $c_{q, m}$ algebraically twice.
Then, $f(c)$ is not null-homotopic in $V -\mathrm{int}V'  \cong T^2 \times I$.
On the other hand,
since $c$ bounds a disk in $S^1 \times I$
whose image under $f$ is contained in $V -\mathrm{int}V'$,
$f(c)$ is null-homotopic in $V -\mathrm{int}V'$.
This is a contradiction.
\QED{Proposition~\ref{hyp ann pair1}}

\begin{figure}[hbt]
\begin{center}
\includegraphics[width=0.9\linewidth]{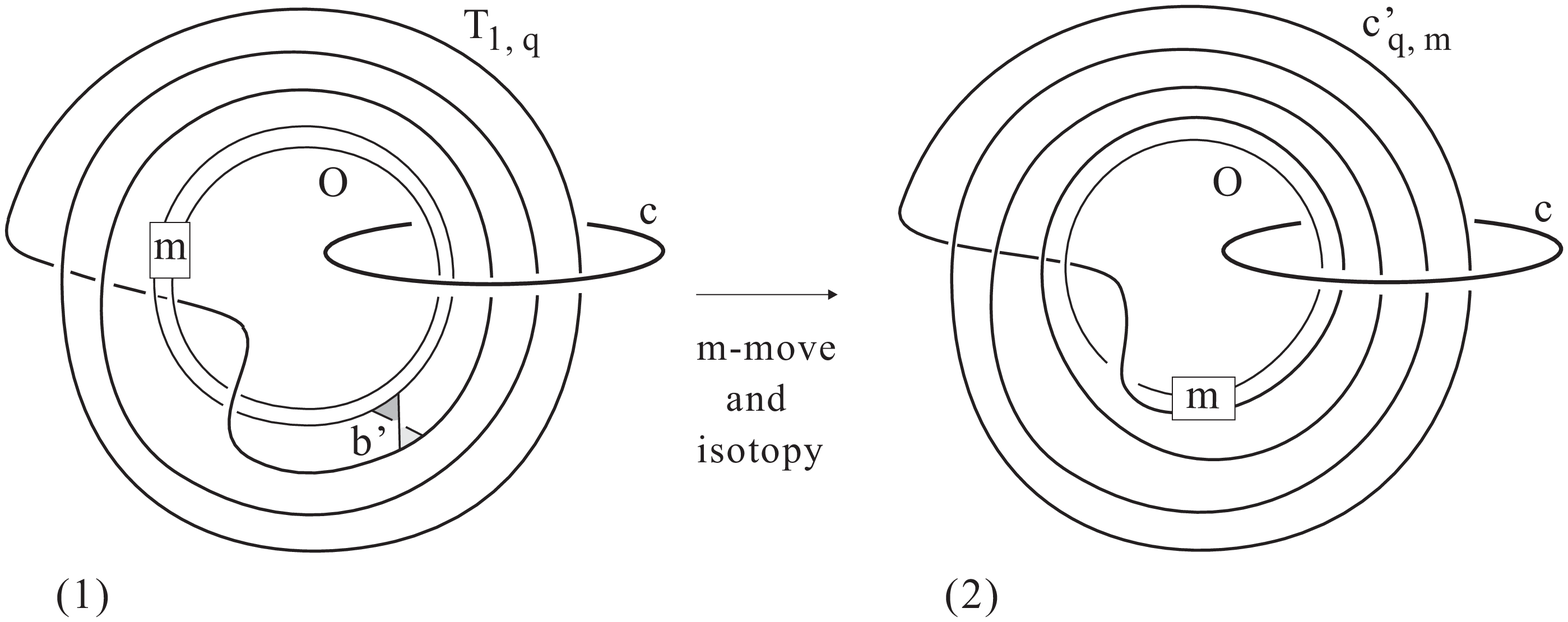}
\caption{Annular pair of seiferters $\{c, c'_{q, m}\}$;\ $q= 3$}
\label{Omhighlinking4}
\end{center}
\end{figure}

\begin{proposition}
\label{hyp ann pair2}
Let $c \cup c'_{q, m}$ be the link obtained from $c \cup T_{1, q}$ in $S^3 -N(O)$ 
by an $m$-move using the band $b'$ in Figure~\ref{Omhighlinking4}(1) 
and an isotopy. 
Assume that $q \ge 1$, $m \ne -1, 0$, and $(m, q) \ne (1, 1)$. 
Then $\{ c, c'_{q, m} \}$ is a hyperbolic annular pair of seiferters for $(O, m)$ 
with $\mathrm{lk}(c, c'_{q, m}) = q+1$. 
\end{proposition}

\begin{figure}[hbt]
\begin{center}
\includegraphics[width=0.9\linewidth]{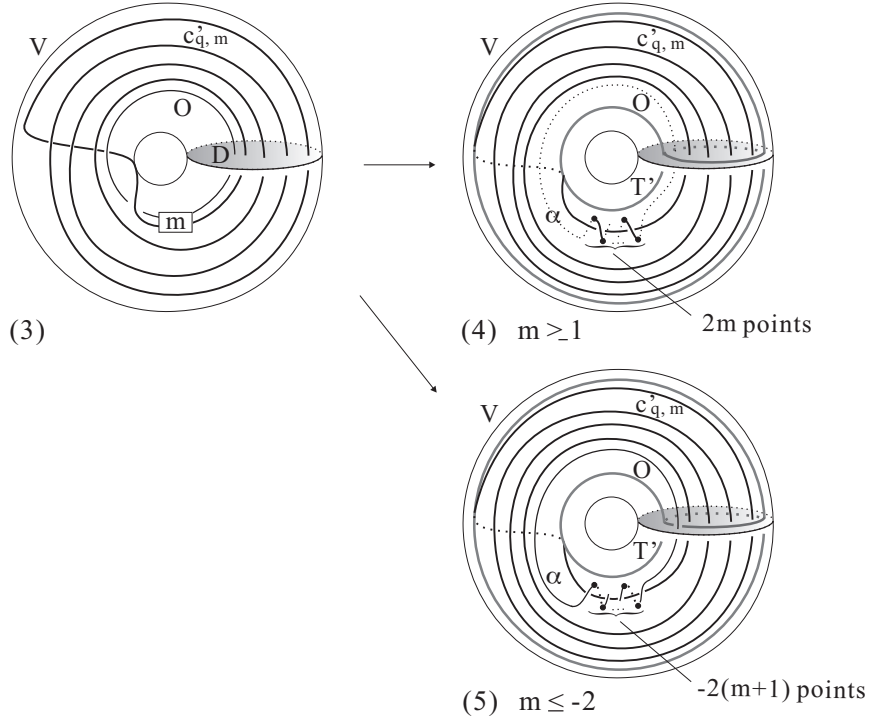}
\caption{In (4), (5),
the intersection points between $O$ and $T'$ are indicated  by ``dots''.}
\label{Omhighlinking5}
\end{center}
\end{figure}

\noindent
\textit{Proof of Proposition $\ref{hyp ann pair2}$}.
Apply the same argument as in the proof of Proposition~\ref{hyp ann pair1} 
with replacement of Claims~\ref{arcs} and \ref{q=3} by 
Claims~\ref{arcs2} and \ref{q=1} below. 
\QED{Proposition~\ref{hyp ann pair2}}

\begin{claim}
\label{arcs2}
\begin{enumerate}
\item
In $V'$ $($resp.\  $V -\mathrm{int}V'$$)$,
the arc $\alpha$ in Figure~\ref{Omhighlinking5}(4) $($resp.\ (5)$)$
is isotopic with $\partial \alpha$ fixed to an arc in $T'$
intersecting $c'_{q, m}$ algebraically twice.
\item
In $V'$ $($resp.\  $V -\mathrm{int}V'$$)$,
each component $\beta$ of $O \cap V'$
$($resp.\ $O \cap (V -\mathrm{int}V')$$)$ other than $\alpha$ is isotopic
with $\partial \beta$ fixed to an arc in $T'$ intersecting
$c'_{q, m}$ once.
\end{enumerate}
\end{claim}

\begin{claim}
\label{q=1}
It holds that $q = 1$, $m \ge 1$, and $f(c_2)$ is the arc
$\alpha ( \subset V' )$ in Figure~\ref{Omhighlinking5}$(4)$.
\end{claim}

\begin{remark}
\label{examples}
Assume that $p \ge 2$, $m \ne 0, \pm 1$.  
Then $\{ c, c_{p+1, m} \}$ in Proposition~\ref{hyp ann pair1} and 
$\{ c, c'_{p-1, m} \}$ in Proposition~\ref{hyp ann pair2} are both 
hyperbolic annular pairs of seiferters for $(O, m)$ with 
$\mathrm{lk}(c, c_{p+1, m}) = \mathrm{lk}(c, c'_{p-1, m}) = p$. 
Since $\{ |\mathrm{lk}(c, O)|,\, |\mathrm{lk}(c_{p+1, m}, O)| \} 
= \{ 1, |1-m| \}$ 
does not coincide with 
$\{ |\mathrm{lk}(c, O)|,\, |\mathrm{lk}(c'_{p-1, m}, O)| \}
= \{ 1, |1+m| \}$,
$\{ c, c_{p+1, m} \}$ and $\{ c, c'_{p-1, m} \}$ are
distinct, annular pairs for $(O, m)$. 
\end{remark}

\section{Seiferters and Hopf pairs for $(T_{p, 2}, m)$}
\label{Tp2m}
\begin{theorem}
\label{Tp2}
For nontrivial torus knots $T_{p, 2}$ $(|p|\ge 3)$,
the following hold.
\begin{enumerate}
\item
Each Seifert surgery $(T_{p, 2}, m)$ 
has a hyperbolic Hopf pair of seiferters. 
\item
A Seifert surgery $(T_{p, 2}, m)$ has a hyperbolic seiferter if 
$m \ne 2p \pm 1$ and $(m, p) \ne (4, 3), (-4, -3)$. 
\end{enumerate}
\end{theorem}

\noindent
\textit{Proof of Theorem~\ref{Tp2}.}
Theorem~\ref{Tp2}(1) follows from
Proposition~\ref{hypseiferterT}(1) below.

Theorem~\ref{Tp2}(2) follows from
Proposition~\ref{hypseiferterT}(2)
if $|p|\ge 5$ and $m \ne 2p$.
The case when $m =2p$ follows from 
the fact that $(T_{p, q}, pq)$ has a hyperbolic seiferter for any nontrivial torus knot $T_{p, q}$
(Claim~\ref{seiferter c+-} and \cite[Lemma~9.1]{MM1}).
The remaining case is when $|p| =3$.
For trefoil knots, various seiferters and annular pairs
are found in \cite{DMMtrefoil}.
For example,
we see from Remark~\ref{T325}(1) that
$(T_{3, 2}, m)$ (resp.\ $(T_{-3, 2}, m)$) has a hyperbolic seiferter
if $m \ne 4$ (resp.\ $m \ne -4$). 
This shows Theorem~\ref{Tp2}(2) with $|p| =3$.
\QED{Theorem~\ref{Tp2}}

\begin{figure}[hbt]
\begin{center}
\includegraphics[width=0.4\linewidth]{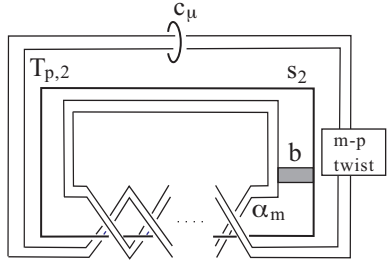}
\caption{$c_m = s_2\,\natural_b\,\alpha_m$}
\label{seiferterTp2}
\end{center}
\end{figure}

\begin{proposition}
\label{hypseiferterT}
Let $c_m$ be the knot obtained from the basic seiferter $s_2$
for $(T_{p, 2}, m)$ $(|p| \ge 3)$ by an $m$-move using the band $b$
described in Figure~\ref{seiferterTp2}.
Then the following hold.
\begin{enumerate}
\item
For the meridional seiferter $c_{\mu}$ as in Figure~\ref{seiferterTp2},
$\{ c_{\mu}, c_m \}$ is a hyperbolic Hopf pair of seiferters
for $(T_{p, 2}, m)$.
\item 
The knot $c_m$ is a hyperbolic seiferter for $(T_{p, 2}, m)$
if $|p| \ge 5$ and $m \ne 2p, 2p \pm 1$.
\end{enumerate}
\end{proposition}

\begin{remark}
\label{m=2p,2p+1}
If $m =2p$,
then $c_m$ in Proposition~\ref{hypseiferterT}
is the same as the basic seiferter $s_2$ for $T_{p, 2}$.
If $m =2p \pm 1$, then
$c_m$ is $(1, \frac{p \pm 1}{2})$ cable of $s_p$ for $T_{p, 2}$,
i.e.\ $c_m$ is the seiferter $s_{p, \pm 1}$ 
for $(T_{p, 2}, 2p \pm 1)$ defined in
\cite[Corollary~3.15(2)]{DMM1}. 
\end{remark}

\noindent
\textit{Proof of Proposition~\ref{hypseiferterT}.}
In (1) we may assume that $p \ge 3$
because the corresponding result for $p \le -3$ can be derived 
by taking mirror images. 
For the same reason we may assume $p \ge 5$ in (2).

(1) The sequence of isotopies
in Figures~\ref{isotopyTp2cmucm} and \ref{isotopyTp2cmucm_2}
shows that $c_m$ is a trivial knot. 
Since $c_m$ is obtained from $s_2$ by an $m$--move and 
$T_{p, 2}$ is a nontrivial knot, 
$c_m$ is a seiferter for $(T_{p, 2}, m)$
by Proposition~2.19(3) in \cite{DMM1}.  
Furthermore,  
since $\{ c_{\mu}, s_2 \}$ is a pair of seiferters for 
$(T_{p, 2}, m)$ and the band $b$ is disjoint from $c_{\mu}$
in Figure~\ref{seiferterTp2}, 
$\{c_{\mu}, c_m \}$ is a pair of seiferters 
\cite[Lemma~2.25(2)]{DMM1}. 
The last figure of Figure~\ref{isotopyTp2cmucm_3} shows
that $\{c_{\mu}, c_m \}$ is a Hopf pair of seiferters. 
Let us verify that no annulus cobounded by $c_{\mu}$ and $c_m$
intersects $T_{p, 2}$ if $m \ne p \pm 1$.
This implies that
$\{c_{\mu}, c_m \}$ $(m \ne p \pm 1)$ is not irrelevant
and thus an annular pair (Remark~\ref{rem:irrelevant});
in particular, $\{c_{\mu}, c_0\}$ is an annular pair of seiferters.
Since $| \mathrm{lk}(c_m, T_{p, 2}) | = |m-p| \ne 1 
= | \mathrm{lk}(c_{\mu}, T_{p, 2}) |$,
$c_m$ is not homologous to $c_{\mu}$ in $S^3 -T_{p, 2}$.
It follows that $c_{\mu}$ and $c_m$ does not cobound an annulus
disjoint from $T_{p, 2}$, as desired.

\begin{figure}[hbt]
\begin{center}
\includegraphics[width=1.0\linewidth]{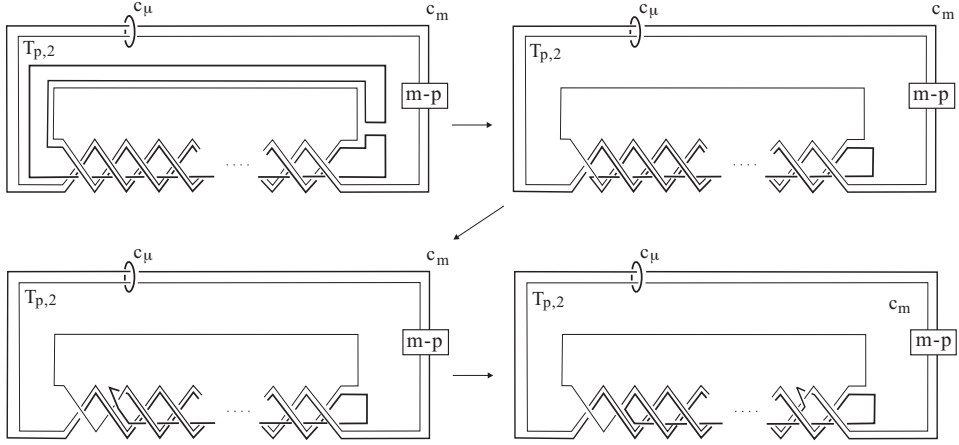}
\caption{Isotopy of $T_{p, 2} \cup c_{\mu} \cup c_m$}
\label{isotopyTp2cmucm}
\end{center}
\end{figure}

\begin{figure}[hbt]
\begin{center}
\includegraphics[width=1.0\linewidth]{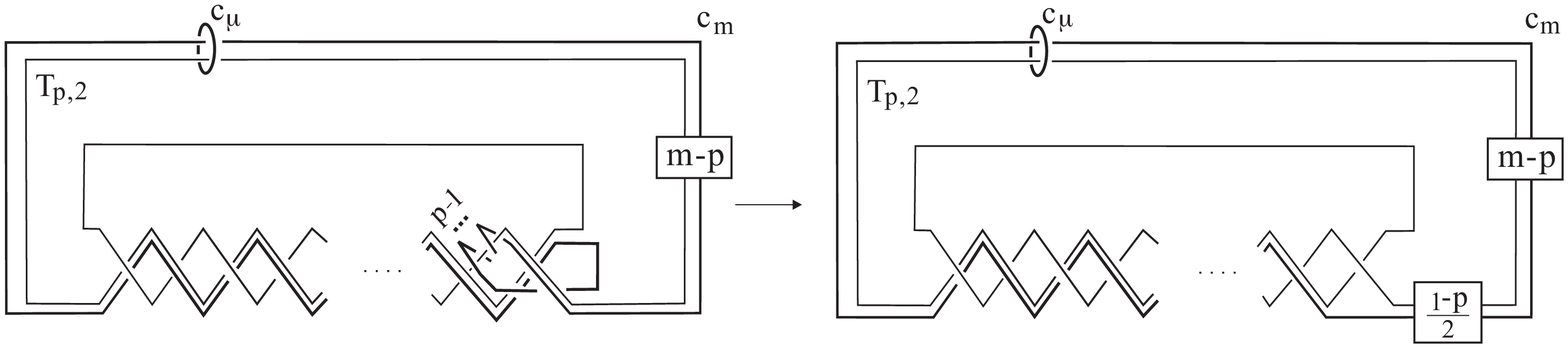}
\caption{Continued from Figure~\ref{isotopyTp2cmucm}}
\label{isotopyTp2cmucm_2}
\end{center}
\end{figure}

\begin{figure}[hbt]
\begin{center}
\includegraphics[width=1.0\linewidth]{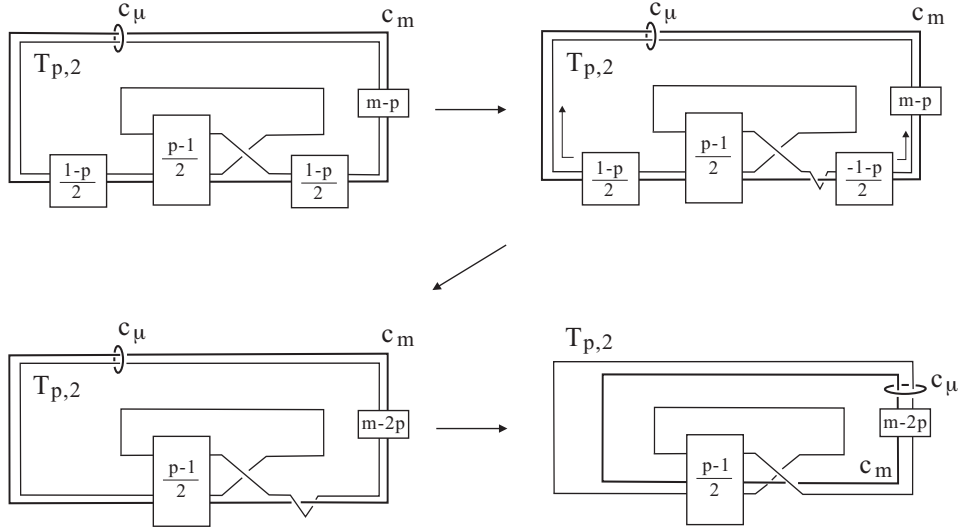}
\caption{Continued from Figure~\ref{isotopyTp2cmucm_2}}
\label{isotopyTp2cmucm_3}
\end{center}
\end{figure}

Let us show that $\{c_{\mu}, c_m\}$ is a hyperbolic annular pair
for any $m$.
In the last figure of Figure~\ref{isotopyTp2cmucm_3},
$(-m)$--twist along $c_{\mu}$ changes $T_{p, 2} \cup c_m$
to $T_{p, 2} \cup c_0$.
Hence, there is an orientation preserving homeomorphism 
from $S^3 - T_{p, 2} \cup c_{\mu} \cup c_m$ to  
$S^3 - T_{p, 2} \cup c_{\mu} \cup c_0$. 
Thus it is sufficient to show that
$S^3 - T_{p,2} \cup c_{\mu} \cup c_0$ is hyperbolic. 
Since $c_{\mu} \cup c_0$ is isotopic to $c_{\mu} \cup s_2$
in $T_{p, 2}(0)$,
$c_{\mu}$ and $c_0$ are exceptional fibers
of indices $2p$ and 2, respectively
in the small Seifert fiber space $T_{p, 2}(0)$
over $S^2(2, p, 2p)$.
Then apply Theorem~3.24 in \cite{DMM1} to
the annular pair $\{c_{\mu}, c_0\}$.
We see that
$\{c_{\mu}, c_0\}$ is a hyperbolic annular pair or 
a basic annular pair,
i.e.\ a pair of basic seiferters $c_{\mu}, s_2, s_p$
as drawn in Figure~\ref{basicseiferters}. 
However, the latter does not occur
because $|\mathrm{lk}(c_{\mu}, c_0)| = 1$.
Hence,  $\{c_{\mu}, c_0\}$ and thus $\{c_{\mu}, c_m\}$ are hyperbolic annular pairs for $(T_{p, 2}, m)$. 
\par

(2) Assume that $m\ne 2p, 2p\pm 1$, and $p \ge 5$. 
As shown in (1) $c_m$ is
isotopic to $s_2$ in $T_{p, 2}(m)$,
and thus an exceptional fiber of index $2$ in $T_{p, 2}(m)$,
a Seifert fiber space over $S^2(2, p, |2p -m|)$.
Then, \cite[Corollary~3.15]{DMM1} 
shows that
$c_m$ is either a basic or a hyperbolic seiferter.
Therefore, Claim~\ref{cm not basic} below implies that
$c_m$ is a hyperbolic seiferter for $(T_{p, 2}, m)$,
as claimed in Proposition~\ref{hypseiferterT}(2).
\QED{Proposition~\ref{hypseiferterT}}

\begin{claim}
\label{cm not basic}
The seiferter $c_m$ in Figure~\ref{seiferterTp2}
is not a basic seiferter for $(T_{p, 2}, m)$.
\end{claim}

\noindent
\textit{Proof of Claim~\ref{cm not basic}.}
We observe the following
from the last figure in Figure~\ref{isotopyTp2cmucm_3}.
{
\renewcommand{\labelenumi}{ (\roman{enumi}) }
\begin{enumerate}
\item
The seiferter $c_{2p}$ is the same as
the basic seiferter $s_2$ for $T_{p, 2}$.
\item
The link $T_{p, 2} \cup c_m$ is obtained from $T_{p, 2} \cup s_2$ 
after $(m -2p)$--twist along $c_{\mu}$.
\end{enumerate}
}

Let $M = S^3 -\mathrm{int}N(T_{p, 2} \cup c_{\mu} \cup c_{2p})$;
$M$ is proved to be hyperbolic
in the proof of Proposition~\ref{hypseiferterT}(1).
We see from observations~(i), (ii) above that
the $\frac{1}{2p -m}$--Dehn filling $M( \frac{1}{2p -m} )$
along $\partial N( c_{\mu} )$ is homeomorphic to
$S^3 -\mathrm{int}N(T_{p, 2} \cup c_m)$, and
$M(\frac{1}{0}) \cong S^3 -\mathrm{int}N(T_{p, 2} \cup s_2)$
is a Seifert fiber space.
Now assume for a contradiction
that $c_m$ is a basic seiferter for $T_{p, 2}$;
then $S^3 -\mathrm{int}N(T_{p, 2}  \cup c_m)$ is Seifert fibered.
By \cite[Corollary~1.2]{GW} we obtain $|2p -m | \le 3$.
Since $m \ne 2p, 2p \pm 1$, it follows $|2p -m | =2$ or $3$. \par

Assume $|2p -m| =2$; then $|\mathrm{lk}(c_m, T_{p, 2})|
= m -p = p +2(>0)$ or $p -2(>0)$.
If $c_m$ is the same as $s_2$, then
we have $|\mathrm{lk}(c_m, T_{p, 2})|
= |\mathrm{lk}(s_2, T_{p, 2})| =p$, a contradiction.
If $c_m$ is the same as $s_p$, then we have
$|\mathrm{lk}(c_m, T_{p, 2})| =2$,
so that $p =0$ or $4$. 
This is not the case 
because $p$ is an odd integer. 
If $c_m$ is the same as $c_{\mu}$,
then since $|\mathrm{lk}(c_{\mu}, T_{p, 2})| =1$,
we obtain $p =-1, 3$.
This contradicts the assumption $p \ge 5$. 

Assume $|2p -m| =3$;
then $|\mathrm{lk}(c_m, T_{p, 2})| = m -p = p +3(>0)$ or $p -3(>0)$.
By comparing linking numbers as above,
we can see that $c_m$ is distinct from $s_2, c_{\mu}$,
and thus $c_m$ is the same as $s_p$ and $p=5$.
Since $c_m$ is also isotopic to $s_2$ in $T_{p, 2}(m)$,
$T_{p, 2}(m) -\mathrm{int}N(s_2)$
(a Seifert fiber space over $D^2(p, |2p - m|) = D^2(5, 3)$) is homeomorphic to
$T_{5, 2}(m) -\mathrm{int}N(s_5)$
(a Seifert fiber space over $D^2(2, 3)$). 
This homeomorphism does not preserve Seifert fibrations
up to isotopy,
a contradiction to \cite[Theorem~VI.18]{J}.
\QED{Claim~\ref{cm not basic}}

As for $T_{-3, 2}$
we find various seiferters and annular pairs of seiferters in 
\cite{DMMtrefoil}. 

\begin{figure}[hbt]
\begin{center}
\includegraphics[width=0.32\linewidth]{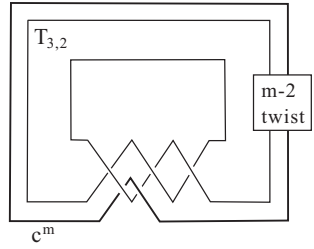}
\caption{Seiferter $c^m (= c_{m +2})$ for $(T_{3, 2}, m)$}
\label{seiferterT32}
\end{center}
\end{figure}

\begin{proposition}[\textbf{\cite{DMMtrefoil}}]
\label{T32}
Take the knot $c^m$ in $S^3 -T_{3, 2}$ illustrated
in Figure~\ref{seiferterT32};
then $c^m$ is a hyperbolic seiferter 
for $(T_{3, 2}, m),\ (T_{3, 2}, m+1)$, and $(T_{3,2}, m+2)$ 
except when $m = 2, 3, 4, 5$. 
In particular, 
$(T_{3, 2}, m)$ has a hyperbolic seiferter if $m \ne 4, 5$. 
\end{proposition}

\begin{remark}
\label{T325}
\begin{enumerate}
\item By setting $n =2$ in Proposition~\ref{non lens-equivalent2}, 
we see that $(T_{3, 2}, 5)$ has a hyperbolic seiferter.
This together with Proposition~\ref{T32} shows that
$(T_{3, 2}, m)$ has a hyperbolic seiferter
for $m \ne 4$.
\item
The seiferter $c^m$ in Figure~\ref{seiferterT32} 
for $(T_{3, 2}, m)$ is isotopic in $S^3 -T_{3, 2}$
to the seiferter
$c_{m +2}$ for $(T_{3, 2}, m+2)$ in Figure~\ref{seiferterTp2}.
\end{enumerate}
\end{remark}

\section{Seiferters not originating in Seifert fibrations of torus knot spaces} 
\label{section:m-equivalence classes}

As shown in Proposition~\ref{seiferter for T}, 
if $m \ne pq, pq \pm 1$ and $T_{p, q}(m)$ is not a prism manifold, 
then any seiferter for $(T_{p, q}, m)$ is $m$--equivalent to a basic seiferter or a regular fiber of 
$S^3 - N(T_{p, q})$.  
On the contrary, as shown in this section,
there exist seiferters for $(T_{p, q}, m)$
which cannot be obtained from basic seiferters or regular fibers
by a sequence of $m$--moves.
In fact, for all $T_{p, q}$ but $T_{\pm 3, 2}$
the degenerate Seifert surgery
$(T_{p, q}, pq)$ has a hyperbolic seiferter not
$pq$--equivalent to a basic seiferter or a regular fiber
in $S^3 -N(T_{p, q})$;
for some $T_{p, q}$ Seifert surgeries $(T_{p, q}, m)$ where
$m =pq +1$ or $pq -1$ have such seiferters.
Examples of the former statement will be given in
Proposition~\ref{non degenerate-equivalent},
and those of the latter
in Propositions~\ref{non lens-equivalent},
\ref{non lens-equivalent2}.

\begin{proposition}
\label{non degenerate-equivalent}
Each Seifert surgery 
$(T_{p, q}, pq)$ $(|p| > q \ge 2)$
where $(p, q) \ne (\pm 3, 2)$ has
a hyperbolic seiferter which is not
$pq$--equivalent to any basic seiferter for $T_{p, q}$ or
a regular fiber of $S^3 - N(T_{p, q})$. 
Furthermore, 
if $|p+q|$ and $|p-q|$ are both greater than one, 
then $(T_{p, q}, pq)$ has at least two such hyperbolic seiferters. 
\end{proposition}

\begin{figure}[htbp]
\begin{center}
\includegraphics[width=0.75\linewidth]{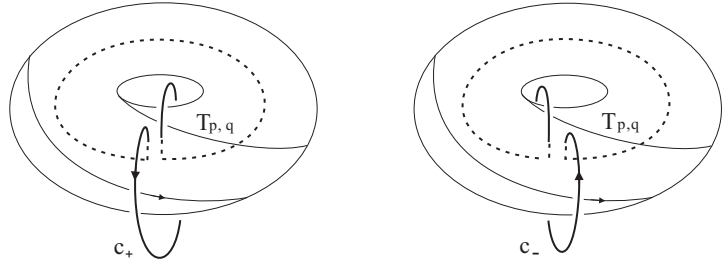}
\caption{Hyperbolic seiferters $c_{+}$ and $c_{-}$ for $(T_{p, q}, pq)$}
\label{pq}
\end{center}
\end{figure}

\noindent
\textit{Proof of Proposition~\ref{non degenerate-equivalent}.}
Let $c_{+}, c_{-}$ be the knots in the exterior of a nontrivial torus knot
$T_{p, q}$ as described in Figure~\ref{pq}. 
The link $T_{p, q} \cup c_{+}$ is exactly the same as the link
$T_{p, q} \cup c$ in \cite[Figure~4.2]{DMM1}; 
see also \cite[Fig.\ 13]{MM1}. 
Note that $\mathrm{lk}(c_{+}, T_{p, q}) = p +q$ and 
$\mathrm{lk}(c_{-}, T_{p, q}) = p - q$. 
The result on $c_{+}$ in Claim~\ref{seiferter c+-} below
is essentially obtained in  \cite[Lemma~9.1]{MM1}. 
Since the link $T_{p, q} \cup c_{-}$ is the mirror image of
$T_{-p, q} \cup c_{+}$, the statement on $c_{-}$ also holds.

\begin{claim}
\label{seiferter c+-}
The knots $c_{\pm}$ are seiferters for $(T_{p, q}, pq)$.
Each of $c_{\pm}$ is a degenerate Seifert fiber in
$T_{p, q}(pq)$ such that $T_{p, q}(pq) -\mathrm{int}N(c_{\pm})$
is a Seifert fiber space over the disk with two exceptional fibers of indices $|p|, q$.
Furthermore, if $|p +q| \ne 1$ $($resp.\  $|p -q| \ne 1$$)$, 
then $c_{+}$ $($resp.\  $c_{-}$$)$ is a hyperbolic seiferter for $(T_{p, q}, pq)$; 
otherwise, $c_{+}$ $($resp.\  $c_{-}$$)$ is a meridian of $T_{p, q}$.
\end{claim}

Since there are no $p, q$ $(|p| > q \ge 2)$ satisfying
$|p +q| =|p -q| =1$,
at least one of $c_{+}, c_{-}$ is a hyperbolic seiferter for $(T_{p, q}, pq)$. 
Set $c = c_+$ if $|p +q| \ne 1$, and otherwise $c = c_-$. 

Let us show that
$c$ is not $pq$--equivalent to any basic seiferter 
if $(p, q) \ne (3, \pm 2)$.
If $c$ were $pq$--equivalent to $s_p$ (resp.\  $s_q$),
then the Seifert fiber space $T_{p, q}(pq) - \mathrm{int}N(c)$
would be homeomorphic to 
$T_{p, q}(pq) - \mathrm{int}N(s_p)
\cong S^1\!\times\! D^2 \sharp L(q, p)$ 
(resp.\ 
$T_{p, q}(pq) - \mathrm{int}N(s_q)
\cong S^1\!\times\! D^2 \sharp L(p, q)$), 
a contradiction to Claim~\ref{seiferter c+-}.
If $c$ were $pq$--equivalent to a meridional seiferter $c_{\mu}$, 
then \cite[Proposition~2.22(1)]{DMM1}  
would show that 
$\mathrm{lk}(c, T_{p, q}) = \pm 1 + xpq$ for some integer $x$.  
Since $\mathrm{lk}(c, T_{p, q}) = p \pm q$, 
a simple computation shows $(p, q) = (\pm 3, 2)$,
a contradiction to our assumption. 
If $c$ were $pq$--equivalent to a regular fiber $t$ in 
$S^3 - N(T_{p, q})$, 
then the Seifert fiber space $T_{p, q}(pq) - \mathrm{int}N(c)$
would be homeomorphic to 
$T_{p, q}(pq) - \mathrm{int}N(t)
\cong S^1\!\times\! D^2 \sharp L(p, q) \sharp L(q, p)$, 
a contradiction.  

Suppose that 
$|p+q|$ and $|p-q|$ are both greater than one. 
We then see that  
$c_+$ and $c_-$ are both hyperbolic seiferters for $(T_{p, q}, pq)$ 
with the required property. 
Since $|\mathrm{lk}(c_+, T_{p, q})| = |p+q| \ne |p-q| = |\mathrm{lk}(c_-, T_{p, q})|$, 
$c_+$ and $c_-$ are distinct seiferters. 
\QED{Proposition~\ref{non degenerate-equivalent}}

\begin{remark}
\label{c+-}
$(T_{p, q}, pq)$ may have a hyperbolic seiferter other than $c_{+}$ and $c_{-}$.
For example, $(T_{3, 5}, 15)$ has a hyperbolic seiferter $c$ such that 
$\mathrm{lk}(c, T_{3, 5}) = 4$.
Since $4 \ne |3 \pm 5|$, $c$ is neither $c_{+}$ nor $c_{-}$.
See \cite[Remark~9.20(1)]{DMM1}.  
\end{remark}

A seiferter for $(T_{p, q}, pq)$ which is not $pq$--equivalent to any basic seiferter or a regular fiber of 
$S^3 - N(T_{p, q})$ arises because of non-uniqueness of degenerate Seifert fibrations of $T_{p, q}(pq)$. 
Similarly,
non-uniqueness of Seifert fibrations of lens spaces
make it possible for some lens surgeries to have such seiferters.  

\begin{proposition}
\label{non lens-equivalent}
The lens surgery 
$(T_{2n + 1, n}, n(2n + 1) -1)$ $(n \ge 2)$ has a hyperbolic seiferter 
which is not $(n(2n + 1) -1)$--equivalent to 
any basic seiferter for $T_{2n + 1, n}$ or a regular fiber of 
$S^3 - N(T_{2n + 1, n})$. 
\end{proposition}

\begin{figure}[htb]
\begin{center}
\includegraphics[ width=1.0\linewidth]{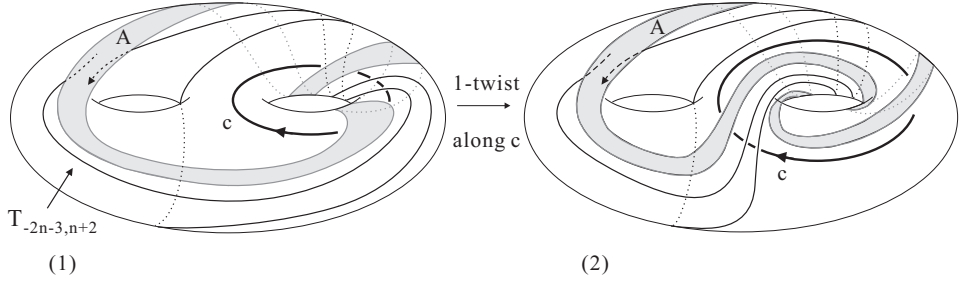}
\caption{An $n$--Dehn twist is performed along the annulus $A$.}
\label{seiferterT2n+1ntwist}
\end{center}
\end{figure}

\begin{figure}[htb]
\begin{center}
\includegraphics[ width=1.0\linewidth]{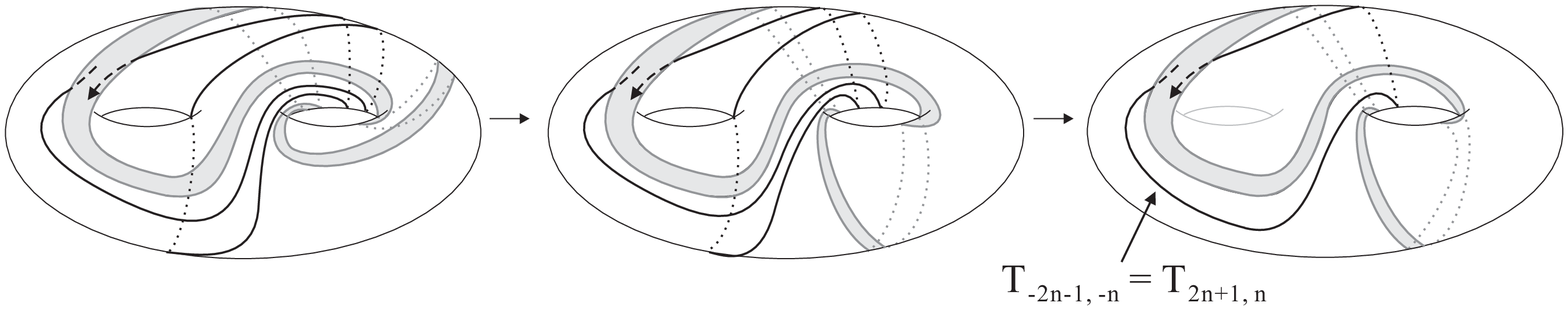}
\caption{}
\label{seiferterT2n+1nisotopy}
\end{center}
\end{figure}

\noindent
\textit{Proof of Proposition~\ref{non lens-equivalent}.}
In \cite[Proposition~3.7]{DMM2}, 
we prove that $c$ described in 
Figure~\ref{seiferterT2n+1ntwist}(1) is a seiferter for the lens surgery 
$(T_{-2n-3, n+2}, (-2n-3)(n+2)+1)$. 
Twisting $T_{-2n-3, n+2}$ once along the seiferter $c$, 
we obtain Figure~\ref{seiferterT2n+1ntwist}(2). 
Figure~\ref{seiferterT2n+1nisotopy} demonstrates that the image of 
$T_{-2n-3, n+2}$ after the twisting is $T_{2n+1, n}$; 
since $\mathrm{lk}(c, T_{-2n-3, n+2}) = 2n+2$, 
the resulting surgery slope is $(-2n-3)(n+2)+1 + (2n+2)^2 = n(2n+1) -1$. 
Thus we obtain the lens surgery $(T_{2n+1, n}, n(2n+1) -1)$ for which $c$ remains a seiferter. 
Note that $ \mathrm{lk}(c, T_{2n+1, n}) = \mathrm{lk}(c, T_{-2n-3, n+2}) = 2n+2$.

Let $(K_p, m_p)$ be the Seifert surgery obtained from
$(T_{2n + 1, n}, n(2n + 1) -1)$ after $p$--twist along $c$.
Then $(K_p, m_p)$ $(p\in \mathbb{Z})$ are Berge's lens surgeries on
Type~III knots. 
Proposition~3.8 in \cite{DMM2} shows that
each lens space $K_p(m_p)$ 
has a Seifert fibration $\mathcal{F}$ over $S^2$ 
such that $\mathcal{F}$ has two exceptional fibers and
$c$ (the image of $c$ after twisting) is a regular fiber
of $\mathcal{F}$.
Hence,
$K_p(m_p) -\mathrm{int}N(c)$ is a Seifert fiber space over the disk
with two exceptional fibers.  

Since $\mathrm{lk}(c, T_{2n +1, n}) = 2n+2
\not \in \{1,\, n,\, 2n+1 \}$, 
the seiferter $c$ is not a basic seiferter for $T_{2n + 1, n}$.
Then, if $c$ were not a hyperbolic seiferter for
$(K_0, m_0) = (T_{2n + 1, n}, n(2n + 1) -1)$,
case~(2), (4), (5), (6), or (7) in
Corollary~3.15 in \cite{DMM1} would occur.

In these cases, $c$ is a $(1, x)$ cable $(|x| \ge 2)$
of an unknotted solid torus
$V$ in $S^3$,
$K_0$ is a knot in $U = S^3 -\mathrm{int}V$,
and a Seifert fibration of $K_0(m_0)$ restricts to that of $V$
with $c$ a regular fiber.
Now for a knot k in a 3-manifold $X (\subset S^3)$,
let us denote by $X(k; \gamma)$
the manifold obtained from $X$ after $\gamma$--surgery on $k$.
If $|p| \ge 2$, then $V(c; -\frac{1}{p})$ has a Seifert fibration
over the disk in which a core of the filled solid torus
is an exceptional fiber of index $|px +1|$
and a core of $V$ is another exceptional fiber of index $|x|$.
In cases~(2), (4), (5), and (7),
$U(K_0; m_0)$ has a Seifert fibration over the disk with at most
two exceptional fibers.
Hence $K_p(m_p) = U(K_0; m_0) \cup V(c; -\frac{1}{p})$ is either
a Seifert fiber space with more than two exceptional fibers
or a lens space which has
two exceptional fibers with $c$ (the image of $c$ after $p$--twist)
one of them.
The former case contradicts the fact that
$(K_p, m_p)$ is a lens surgery for any $p$.
The latter implies $K_p(m_p) -\mathrm{int}N(c)$
is a solid torus, a contradiction.
The remaining case is (6) in 
Corollary~3.15 in \cite{DMM1}.
In this case, $K_0(m_0) -\mathrm{int}N(c)$ is a Seifert fiber space 
over the M\"obius band with one exceptional fiber,
a contradiction.
It follows that $c$ is a hyperbolic seiferter for $(K_0, m_0)$.

Finally we show that 
$c$ is not $m_0$--equivalent to 
any basic seiferter for $K_0 = T_{2n + 1, n}$ or a regular fiber in 
$S^3 - N(K_0)$, where $m_0 = n(2n + 1) -1$. 
If $c$ is $m_0$--equivalent to a basic seiferter 
$s_{2n + 1}$ or $s_{n}$ for $K_0$, 
then $c$ is an exceptional fiber in the lens space $K_0(m_0)$.
It follows that $K_0(m_0) -\mathrm{int}N(c)$ is a solid torus,
a contradiction.
Let us suppose that $c$
is $m_0$--equivalent to $c_{\mu}$.
Then, since $|\mathrm{lk}(c_{\mu}, K_0)| =1$,
\cite[Proposition~2.22(1)]{DMM1} implies 
$\mathrm{lk}(c, K_0) = \pm 1 +xm_0$ for some integer $x$. 
On the other hand, 
$\mathrm{lk}(c, K_0) = 2n+2$. 
We thus have $\pm 1 + xm_0 = 2n+2$, where $n \ge 2$.  
Then $x = \frac{2n+1}{2n^2 +n-1}$ or 
$\frac{2n+3}{2n^2 +n-1}$;   
these cannot be integers because $2n^2 +n-1 > 2n+3 > 2n+1 > 0$ for $n \ge 2$. 
Hence $c$ cannot be $m_0$--equivalent to $c_{\mu}$. 
Let us show that the seiferter $c$ for $(K_0, m_0)$
is not $m_0$--equivalent to 
a regular fiber of $S^3 - N(K_0)$.   
Since the linking number between $T_{2n +1, n}$ and
a regular fiber of $S^3 - N(T_{2n + 1, n})$ is $\pm n(2n +1)$.
We obtain $\pm n(2n+1) + xm_0 = 2n+2$ 
for some integer $x$
(\cite[Proposition~2.22(1)]{DMM1}),
where $n \ge 2$.  
Then $x = -1 + \frac{2n+1}{2n^2 +n-1}$ or 
$1 + \frac{2n+3}{2n^2 +n-1}$, 
which cannot be integers for any $n$.  
Hence $c$ cannot be $m_0$--equivalent to a regular fiber 
in $S^3 - N(K_0)$. 
\QED{Proposition~\ref{non lens-equivalent}}

\begin{proposition}
\label{non lens-equivalent2}
The lens surgery 
$(T_{2n - 1, n}, n(2n - 1) -1)$ $(n \ge 2)$ has a hyperbolic seiferter 
which is not $(n(2n - 1) -1)$--equivalent to 
any basic seiferter for $T_{2n - 1, n}$ or a regular fiber of 
$S^3 - N(T_{2n - 1, n})$. 
\end{proposition}

\begin{figure}[htb]
\begin{center}
\includegraphics[ width=1.0\linewidth]{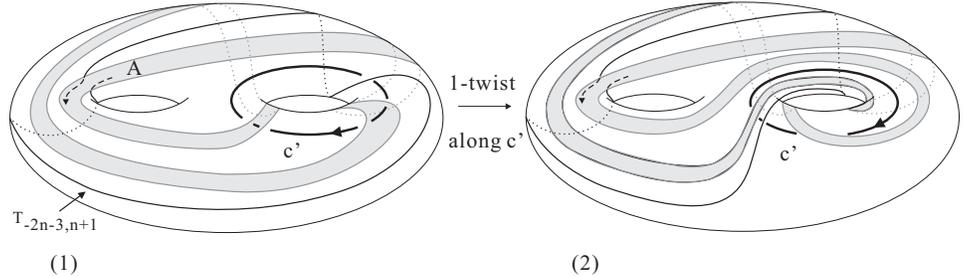}
\caption{An $n$--Dehn twist is performed along the annulus $A$.}
\label{seiferterT2n-1ntwist}
\end{center}
\end{figure}

\begin{figure}[htb]
\begin{center}
\includegraphics[ width=1.0\linewidth]{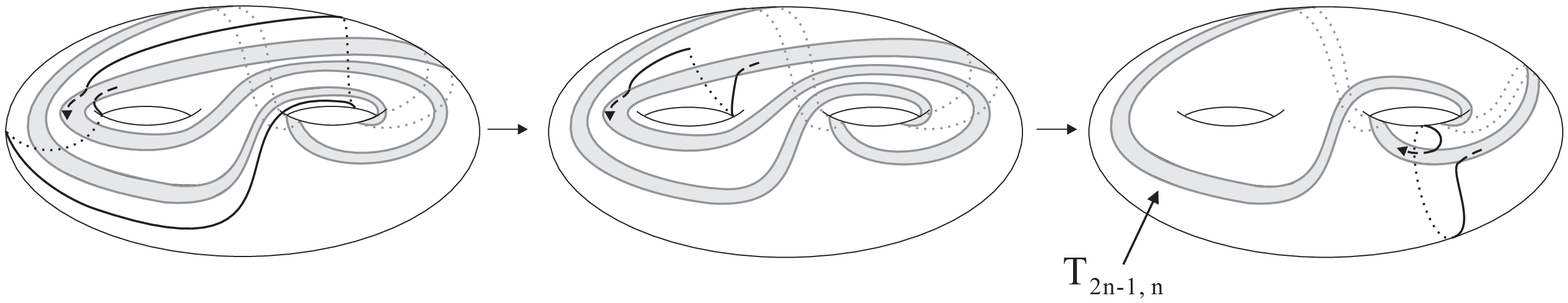}
\caption{}
\label{seiferterT2n-1nisotopy}
\end{center}
\end{figure}

\noindent
\textit{Proof of Proposition~\ref{non lens-equivalent2}.}
In \cite[Section 4]{DMM2}, 
we prove that $c'$ described in 
Figure~\ref{seiferterT2n-1ntwist}(1) is a seiferter for the lens surgery 
$(T_{-2n-3, n+1}, (-2n-3)(n+1)+1)$.  
Twisting $T_{-2n-3, n+1}$ once along $c'$, 
we obtain Figure~\ref{seiferterT2n-1ntwist}(2). 
Figure~\ref{seiferterT2n-1nisotopy} demonstrates that the image of 
$T_{-2n-3, n+1}$ after the twisting is $T_{2n-1, n}$; 
since $\mathrm{lk}(c', T_{-2n-3, n+1}) =  2n+1$, 
the resulting surgery slope is $(-2n-3)(n+1)+1 + (2n+1)^2 = n(2n-1) -1$. 
Thus we obtain a lens surgery $(T_{2n-1, n}, n(2n-1) -1)$ for which $c'$ remains a seiferter. 
Note that $\mathrm{lk}(c', T_{2n-1, n}) = \mathrm{lk}(c', T_{-2n-3, n+1})
=  2n+1$. 

Let $(K_p, m_p)$ be the Seifert surgery obtained from
$(T_{2n-1, n}, n(2n-1) -1)$ after $p$--twist along $c'$.
Then $(K_p, m_p)$ $(p \in \mathbb{Z})$
are Berge's lens surgeries on Type~IV knots. 
In \cite[Section~4]{DMM2} it is shown that each lens space
$K_p(m_p)$ has a Seifert fibration $\mathcal{F}$ over $S^2$ 
such that $\mathcal{F}$ has two exceptional fibers and
$c'$ (the image of $c'$ after twisting) is a regular fiber
of $\mathcal{F}$. 
Since $\mathrm{lk}(c', T_{2n -1, n}) = 2n+1
\not \in \{1,\, n,\, 2n-1 \}$, 
the seiferter $c'$ is not a basic seiferter for $T_{2n - 1, n}$.
Then the argument in the proof of Proposition~\ref{non lens-equivalent} shows that 
$c'$ is a hyperbolic seiferter for the lens surgery $(T_{2n-1, n}, n(2n-1) -1)$, 
and is not $( n(2n-1) -1)$--equivalent to a basic seiferter or a regular fiber of 
$S^3 - N(T_{2n-1, n})$. 
\QED{Proposition~\ref{non lens-equivalent2}}

\begin{remark}
\label{distinct seiferters lens}
The lens surgery $(T_{2n+1, n}, n(2n+1) -1)$ $(n\ge 3)$ has,
other than $c$ in Figure~\ref{seiferterT2n+1ntwist}(2),
a hyperbolic seiferter
which is not $( n(2n+1)-1 )$--equivalent to a basic seiferter or 
a regular fiber of $S^3 - N(T_{2n + 1, n})$.
Let us put $m = n+1$, where  $n \ge 2$; 
then the lens surgery $(T_{-2n-3, n+1}, (-2n-3)(n+1)+1)$
in the proof of Proposition~\ref{non lens-equivalent2}
becomes $(T_{-2m-1, m}, m(-2m-1)+1)$, and 
$\mathrm{lk}(c', T_{-2m-1, m}) = 2m - 1$,
where $c'$ is as in Figure~\ref{seiferterT2n-1ntwist}(1).  
Let $T_{2m+1, m} \cup c'^*$ be the mirror image of
the link $T_{-2m -1, m} \cup c'$.
Writing $n$ for $m (\ge 3)$, 
we have the hyperbolic seiferter $c'^*$ for $(T_{2n+1, n}, n(2n+1)-1)$ which is not $( n(2n +1)-1 )$--equivalent to
a basic seiferter or a regular fiber of 
$S^3 - N(T_{2n +1, n})$.
Since $|\mathrm{lk}(c'^*, T_{2n+1, n})| = 2n-1$
is not equal to $|\mathrm{lk}(c, T_{2n+1, n})| =2n +2$, 
the seiferter $c'^*$ is distinct from $c$. 
\end{remark}

\section{Band sums and seiferters}
\label{single band sum}

For a 2--component link $k_1\  \cup\ k_2$, 
we call a band $b$ connecting $k_1$ and $k_2$
a \textit{trivializing band}
if the band sum $k_1\,\natural_b\,k_2$ is a trivial knot in $S^3$. 
Theorem~\ref{band sum} below determines when we have
a trivializing band 
connecting a torus knot $T_{p, q}$ and its basic seiferters 
$s_p, s_q, c_{\mu}$. 

\begin{theorem}
\label{band sum}
Let $T_{p, q}$ be a nontrivial torus knot with $|p| > q \ge 2$.
Then the following hold.
\begin{enumerate}
\item
There exists a trivializing band connecting $s_q$ and $T_{p, q}$ 
if and only if $q =2$.
\item
There exists a trivializing band connecting $s_p$ and $T_{p, q}$ 
if and only if $(p, q) = (\pm 3, 2)$. 
\item
There exists a trivializing band connecting $c_{\mu}$ and $T_{p, q}$ 
if and only if $(p, q) = (\pm 3, 2)$. 
\end{enumerate}
\end{theorem}

\begin{figure}[htb]
\begin{center}
\includegraphics[width=1.0\linewidth]{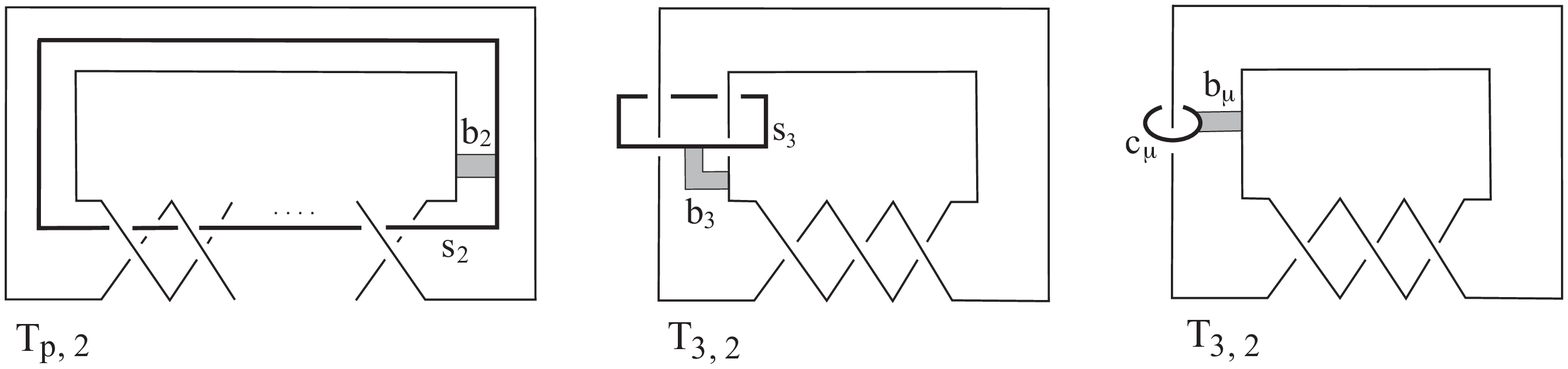}
\caption{Band sums $s_2\,\natural_{b_2}\,T_{p, 2}$, 
$s_3\,\natural_{b_3}\,T_{3, 2}$, and 
$c_{\mu}\,\natural_{b_{\mu}}\,T_{3,2}$ are trivial knots.}
\label{torusband}
\end{center}
\end{figure}

\noindent
\textit{Proof of Theorem~\ref{band sum}.} 
The band sum of $T_{p, 2}$ and $s_2$ described in 
the first figure of Figure~\ref{torusband} is a trivial knot in $S^3$. 
Moreover, if $p = 3$, 
the band sums of $T_{3, 2}$ and basic seiferters $s_3$ and $c_{\mu}$ 
described in the second and the third figures of 
Figure~\ref{torusband} are both trivial knots. 
This fact proves the if parts of Theorem~\ref{band sum}.

The only if part of assertion~(3) is proved in \cite{IM};
it is further shown that
if a band sum of $c_{\mu}$ and $T_{3, 2}$ is a trivial knot,
then the band is isotopic to $b_{\mu}$ in Figure~\ref{torusband}.
Thus it is enough to prove the only if parts of (1), (2).
The proof is done by relating the band sums
to basic seiferters
for the degenerate Seifert surgery $(T_{p, q}, pq)$.

\begin{figure}[htb]
\begin{center}
\includegraphics[width=0.31\linewidth]{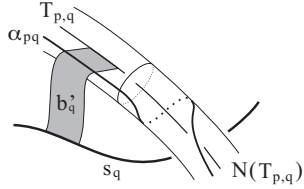}
\caption{Band sum of $s_q$ and $T_{p, q}$, 
and band sum of $s_q$ and $\alpha_{pq}$}
\label{oneband}
\end{center}
\end{figure}

(1) Let $b_q$ be a band connecting $s_q$ and $T_{p, q}$, 
and write $k_q = s_q\,\natural_{b_q}\,T_{p, q}$.
Take a tubular neighborhood of $T_{p, q}$ so that
$N(T_{p, q}) \cap s_q = \emptyset$, and $\partial N(T_{p, q}) \cap b_q$ is an arc.
Let $\alpha_{pq}$ be a simple closed curve on $\partial N(T_{p, q})$
with slope $pq$ such that
$\alpha_{pq} \cap b_q = \partial N(T_{p, q}) \cap b_q$.
Then, $b'_q = b_q -\mathrm{int}N(T_{p, q})$ is a band connecting $\alpha_{pq}$
and $s_q$ (Figure~\ref{oneband}). 
Let $c$ be a knot in $S^3 - N(T_{p, q})$ which is obtained from  
the band sum $s_q\,\natural_{b'_q}\,\alpha_{pq}$ by pushing away from 
$\partial N(T_{p, q})$; 
$c$ is obtained from the basic seiferter $s_q$
by a single $pq$--move using the band $b'_q$. 
Note that $c$ is isotopic to $k_q$ in $S^3$.

Now suppose that the band sum $k_q$ is a trivial knot in $S^3$.
Then, $c$ is a seiferter for $(T_{p, q},\  pq)$; 
moreover, 
since $c$ is isotopic in $T_{p, q}(pq)$ to the basic seiferter $s_q$,
$c$ is a non-degenerate exceptional fiber of index $q$ in $T_{p, q}(pq)$.
Let $V$ be the solid torus $S^3 -\mathrm{int}N(c)$.
We prove the claim below on the position of $T_{p, q}$ in $V$.

\begin{lemma}
\label{position}
The position of $T_{p, q}$ in $V$ is one of the following.
{
\renewcommand{\labelenumi}{ \textup{(\roman{enumi})} }
\begin{enumerate}
\item $T_{p, q}$ is a $(q, p)$ cable of $V$.
\item $T_{p, q}$ is a $(q, p)$ cable of a $(1, s)$ cable of $V$
for some integer $s$ such that $|s| \ge 2$ and $q =sp \pm 1$.
\end{enumerate}
}
\end{lemma}

\noindent
\textit{Proof of Lemma~\ref{position}.}
Since $c$ is a non-degenerate exceptional fiber in $T_{p, q}(pq)
\cong L(p, q) \sharp L(q, p)$,
we have four possibilities (Corollary 3.21(2) and Theorem~3.19(2)(iii) in \cite{DMM1}). 
\begin{itemize}
\item[(i)] $T_{p, q}$ is a $(q, p)$ cable of $V$.
\item[(i$'$)] $T_{p, q}$ is a $(p, q)$ cable of $V$.
\item[(ii)] $T_{p, q}$ is a $(q, p)$ cable of a $(1, s)$ cable of $V$ for some integer $s$ such that $|s| \ge 2$ and $q =sp \pm 1$.
\item[(ii$'$)] $T_{p, q}$ is a $(p, q)$ cable of a $(1, s)$ cable of $V$ for some integer $s$ such that $|s| \ge 2$ and $p =sq \pm 1$.
\end{itemize}

Since $c$ is isotopic to $s_q$ in $T_{p, q}(pq)$, 
we see $V(T_{p, q}; pq) = T_{p, q}(pq) -\mathrm{int}N(c)
\cong T_{p, q}(pq) -\mathrm{int}N(s_q)$. 
This manifold is homeomorphic to $S^1 \times D^2 \sharp L(p, q)$
because $T_{p, q}$ is the $(q, p)$ cable of the solid torus
$S^3 -\mathrm{int}N(s_q)$.
On the other  hand,
in cases~(i$'$) and (ii$'$), 
$V(T_{p, q}; pq) \cong S^1\!\times\! D^2 \sharp L(q, p)$.
Thus cases~(i$'$) and (ii$'$) do not occur.
\QED{Lemma~\ref{position}}

In case (i) in Lemma~\ref{position}, $|\mathrm{lk}(c, T_{p, q})| = p$.
On the other hand,
since $c$ is obtained from $s_q$ by a single $pq$--move,
we have
$\mathrm{lk}(c, T_{p, q})
= \mathrm{lk}(s_q, T_{p, q}) +\varepsilon pq$, 
where $c, s_q, T_{p, q}$ are oriented adequately
and $\varepsilon \in \{ \pm 1 \}$
(\cite[Proposition~2.22(1)]{DMM1}).
Hence, $|p + \varepsilon pq| = |p|$. 
It follows that $q = 0, \pm 2$.
Since $q \ge 2$, 
we obtain $q = 2$ as claimed in Theorem~\ref{band sum}.
\par
Now let us consider case~(ii) in Lemma~\ref{position}
where $T_{p, q}$ is a $(q, p)$ cable of a $(1, s)$ cable of $V$;
then $\mathrm{lk}(c, T_{p, q}) = \pm ps$. 
It follows that  $|p + \varepsilon pq| = |ps|$ 
and thus $|1 + \varepsilon q| = |s|$.
Combining this equality with $|ps - q| = 1$ in case~(ii),
we obtain the inequalities below.
$$|ps| -1 \le |ps \varepsilon +1| \le
|1 +\varepsilon q| +|ps \varepsilon - \varepsilon q| = |s| + 1$$
It follows $|ps| \le |s| +2$. 
Since $|s| \ge 2$, 
$|p| \le 1+ \frac{2}{|s|} \le 2$. 
This contradicts the assumption $|p| > q \ge 2$. 
Assertion~(1) is thus proved.

(2) Starting with a band sum $k_p = s_p\,\natural_{b_p}\,T_{p, q}$,
we follow the argument in (1) with $p$ and $q$ exchanged.
Then, we obtain the same statement as in cases~(i) and (ii)
in Lemma~\ref{position} with $p$ and $q$ exchanged.
The modified case~(i) then leads to $p = 0, \pm 2$.
However, this is impossible because $|p| > q \ge 2$.
The modified case~(ii) leads to the inequality $|qs| \le |s| +2$,
so that $q \le 1 + \frac{2}{|s|}$.
Then, using the fact $|p| > q \ge 2$ and $|s| \ge 2$,
we see that $q =2$ and $|s| =2$.
Since $|1 +\varepsilon q | = |s|$ holds in case~(ii),
$|1 +\varepsilon p| = |s|$ holds in the modified case~(ii).
We then obtain $p = \pm 3$ and $q =2$, as desired
in assertion~(2).
\QED{Theorem~\ref{band sum}}

Theorem~\ref{band sum} implies the following results on
seiferters obtained by $m$--moves.

\begin{theorem}
\label{restriction1}
Let $(T_{p, q}, m)$ be a Seifert surgery on a torus knot $T_{p, q}$
with $|p| > q > 2$. 
Then, there is no seiferter for $(T_{p, q}, m)$ which is obtained
from a basic seiferter by an $m$--move. 
\end{theorem}

\noindent
\textit{Proof of Theorem~\ref{restriction1}.}
Theorem~\ref{band sum} shows that  
all band sums of $T_{p, q}$ $(|p| > q > 2)$ and basic seiferters
for $T_{p, q}$ are nontrivial knots in $S^3$.
Let $\alpha_m$ be a simple closed curve in $\partial N(T_{p, q})$
with slope $m$.
It follows that all band sums of $\alpha_m$ and
basic seiferters for $T_{p, q}$ are nontrivial knots
because $\alpha_m$ is isotopic in $N(T_{p, q})$ to the core $T_{p, q}$.
Thus an arbitrary knot obtained from each basic seiferter
for $T_{p, q}$ by an $m$--move is not a seiferter. 
\QED{Theorem~\ref{restriction1}}

\begin{theorem}
\label{restriction2}
Let $(T_{p, q}, m)$ be a Seifert surgery on a torus knot $T_{p, q}$
with $|p| > q \ge 2$. 
Suppose that $c$ is a seiferter for $(T_{p, q}, m)$  
which is obtained from a regular fiber of $S^3 - N(T_{p, q})$ 
by an $m$--move. 
Then, $c$ is a $(1, m-pq)$ cable of a meridian of $T_{p, q}$,
and thus a non-hyperbolic seiferter for $(T_{p, q}, m)$. 
\end{theorem}

\noindent
\textit{Proof of Theorem~\ref{restriction2}.}
Suppose that $c$ is a seiferter for $(T_{p, q}, m)$  
which is obtained from a regular fiber $t$ in $S^3 - N(T_{p, q})$ 
by an $m$--move
using a band $b (\subset S^3 -\mathrm{int}N( T_{p, q} ) )$; 
$c$ is a knot in $S^3 -N( T_{p, q} )$ obtained by pushing 
the band sum $t\,\natural_b\,\alpha_m$
away from $\partial N(T_{p, q})$,
where $\alpha_m$ is a simple closed curve on $\partial N( T_{p, q} )$
with slope $m$. 
Our purpose is to show that 
$c$ is a $(1, m-pq)$ cable of a meridian of $T_{p, q}$.

In the Seifert fibration of $S^3 -\mathrm{int}N(T_{p, q})$,
take a regular fiber $\alpha_{pq}$ on $\partial N(T_{p, q})$,
which represents the slope $pq$. 
We may assume that
there is a small annulus $M$ on $\partial N(T_{p, q})$ such that
the core curve of $M$ is a meridian of $T_{p, q}$,
and that
$\alpha_m$ and $\alpha_{pq}$ restrict to the same essential arc
in the annulus $\partial N(T_{p, q}) -\mathrm{int}M$.  \par

\begin{figure}[htb]
\begin{center}
\includegraphics[width=0.38\linewidth]{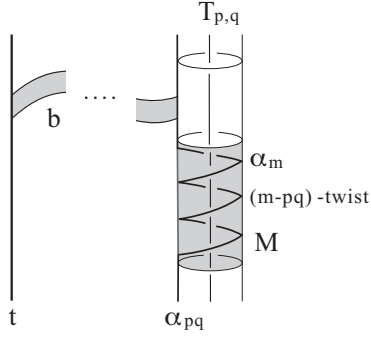}
\caption{A band sum of $\alpha_m$ and a regular fiber $t$}
\label{regularfiberband}
\end{center}
\end{figure}

Now isotope $b$ so that $b \cap \alpha_m$ is contained in 
$\partial N(T_{p, q}) -\mathrm{int}M$,
and take the band sum $t\,\natural_b\,\alpha_{pq}$. 
Note that $t\,\natural_b\,\alpha_m$ and
$t\,\natural_b\,\alpha_{pq}$ coincide outside of $M$
and are isotopic in $S^3$.
Let $c_{pq}$ be a knot obtained by
pushing $t\,\natural_b\,\alpha_{pq}$
away from $\partial N(T_{p, q})$. 
Since $c$ is a trivial knot in $S^3$,  
$c_{pq}$ is also trivial in $S^3$. 
Since $c_{pq}$ is isotopic in $T_{p, q}(pq)$ to the regular fiber $t$, 
$c_{pq}$ is a regular fiber in a degenerate Seifert fibration of
$T_{p, q}(pq)$. 
On the other hand,
Theorem~3.21(1) in \cite{DMM1} shows that
no seiferter for the degenerate Seifert surgery $(T_{p, q}, pq)$
is a regular fiber.
Hence, $c_{pq}$ is not a seiferter.
It follows that $c_{pq}$ is an irrelevant seiferter,
and so bounds a disk in $S^3 -T_{p, q}$
(Remark~\ref{rem:irrelevant}).

On the position of the band $b$ the following holds.

\begin{lemma}
\label{annulus-band}
There exists an annulus $S$ in $S^3 -\mathrm{int}N(T_{p, q})$
such that 
$\partial S = t \cup \alpha_{pq}$ and $b \subset S$. 
\end{lemma}

\begin{figure}[htb]
\begin{center}
\includegraphics[width=0.9\linewidth]{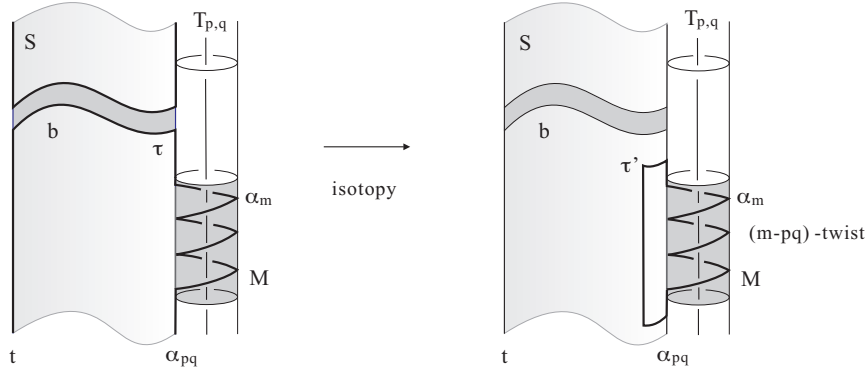}
\caption{$\tau$ is isotopic to $\tau'$.}
\label{regularfiberband2}
\end{center}
\end{figure}

Using the annulus $S$ obtained by this lemma, we complete the proof of
Theorem~\ref{restriction2}.
By an isotopy we may assume further  that 
$S \cap N(T_{p, q}) = \alpha_{pq}$. 
Since $S$ contains the band $b$, 
$t\,\natural_b\,\alpha_m$ is 
the union of the two arcs $\alpha_m \cap M$ and
$\tau = \partial \overline{(S -b)} -\mathrm{int}M$;
$\tau$ is isotopic in $S$ with its end points fixed
to the arc $\tau'$ in Figure~\ref{regularfiberband2}. 
Note that $( \alpha_m \cap M ) \cup \tau'$ is
a $(1, m-pq)$ cable of a meridian of $T_{p, q}$.
This shows that $t\,\natural_b\,\alpha_m$ and thus $c$
is isotopic in $S^3 -T_{p, q}$
to the $(1, m-pq)$ cable of a meridian of $T_{p, q}$, as claimed. 
\QED{Theorem~\ref{restriction2}}

\noindent
\textit{Proof of Lemma~\ref{annulus-band}.} 
Let $A$ be an annulus
in $S^3 - \mathrm{int}N(T_{p, q})$
with $\partial A = t \cup \alpha_{pq}$. 
(Since $t$ and $\alpha_{pq}$ are
regular fibers in the Seifert fibration of
$S^3 - \mathrm{int}N(T_{p, q})$, such an annulus is obtained as
a union of regular fibers.)
Choose orientations of $t$ and $\alpha_{pq}$ which are consistent
in $t\,\natural_b\,\alpha_{pq}$. 
We consider two cases according as 
$t$ and $\alpha_{pq}$ are homologous in $A$ or not. 
First suppose that $t$ and $\alpha_{pq}$ are homologous in $A$. 
Then $\mathrm{lk}(c_{pq}, T_{p, q}) = 
\mathrm{lk}(t, T_{p, q}) + \mathrm{lk}(\alpha_{pq}, T_{p, q})
=2pq \ne 0$. 
However, this contradicts the fact that
$c_{pq}$ bounds a disk in $S^3 -T_{p, q}$.
Now assume that 
$t$ and $\alpha_{pq}$ are not homologous in $A$. 
Then the (adequately oriented) annulus $A$
in $S^3 -\mathrm{int}N( T_{p, q} )$
is a Seifert surface for the oriented link $t \cup \alpha_{pq}$.
Here, 
a Seifert surface $F$ for an oriented link $L$ 
is a compact oriented surface such that
no component of $F$ is closed and $\partial F = L$.
We define $\chi(L)$ to be
the maximal Euler characteristic of all Seifert surfaces for $L$. 
Since $ t \natural_b \alpha_{pq}$ is a trivial knot in $S^3$, 
we see $\chi( t \natural_b \alpha_{pq} ) = 1$.
Since the oriented link $t \cup \alpha_{pq}$ is non-splittable,
it follows $\chi( t \cup \alpha_{pq} ) = \chi( A ) =0$.
Then, the minor revision of \cite[Theorem~1.6]{HiraS} below
shows that
the oriented link $t \cup \alpha_{pq}$ cobounds an annulus $S$ 
containing the band $b$, as claimed.  
By an isotopy we may assume that $S \subset S^3 - \mathrm{int}N(T_{p, q})$. 
\QED{Lemma~\ref{annulus-band}}

\begin{theorem}[\textbf{a minor revision of Theorem~1.6 in \cite{HiraS}}]
\label{theorem1.6}
Let $L$ be an oriented link,
and $b$ a band connecting $($possibly the same$)$ components of $L$
such that
$L$ and $b$ induce opposite orientations to
their intersection $L \cap b$.
Denote the self band sum of $L$ using $b$ by $L_b$,
an oriented link.
Then,
$\chi( L ) \le \chi( L_b ) -1$ if and only if
$L$ has a Seifert surface $S$ such that
$\chi( S ) = \chi( L )$ and $b \subset S$.
\end{theorem}

\noindent
\textit{Proof of Theorem~\ref{theorem1.6}}.
For a Seifert surface $S ( \subset M = S^3 -\mathrm{int}N(L) )$
for $L$, consider the three conditions below.
\begin{enumerate}
\item $S$ is taut in $(M, \partial M)$,
i.e.\ $S$ is incompressible and minimizes the Thurston norm
of $[S, \partial S] \in H_2(M, N)$,
where $N$ is a tubular neighborhood of $\partial S$ in $\partial M$.
\item $\chi( L ) =\chi( S )$
\item $S$ is a minimal genus Seifert surface for $L$,
i.e.\ the sum of the genera of the components of $S$ is minimal. 
\end{enumerate}
Theorem~1.6 in \cite{HiraS} states that
$\chi( L ) \le \chi( L_b ) -1$ if and only if 
$L$ has a minimal genus Seifert surface $S$ such that $b \subset S$.
In the proof
the authors assume that 
$(3) \Rightarrow (2)$ and $(1) \Leftrightarrow (3)$ are true.
However, these are not true;
if a minimal genus Seifert surface $S$ for a link $L$ is
disconnected, then by tubing two components of $S$,
we obtain a minimal genus, compressible Seifert surface $S'$
with $\chi(S') < \chi(L)$.
On the other hand, $(1) \Leftrightarrow (2)$ holds
by \cite[Lemma~1.2]{Sch-Tho}.
By replacing the word ``minimal genus"
in the proof of \cite[Theorem~1.6]{HiraS} with ``taut",
we obtain a proof of Theorem~\ref{theorem1.6}.
\QED{Theorem~\ref{theorem1.6}}

\begin{remark}
\label{algnonsplit}
Among connected Seifert surfaces for a given link, 
a Seifert surface $S$ has minimal genus if and only if 
$\chi(S)$ is maximal. 
Thus,
\cite[Theorem~1.6]{HiraS} holds for
links which have only connected Seifert surfaces.
Theorem~\ref{band sum}(3) is, in fact, proved in \cite{IM} by using
Theorem~1.6 in \cite{HiraS}.
However, in the proof Theorem~1.6 in \cite{HiraS} is applied
only to links with only connected Seifert surfaces.
\end{remark}

\begin{corollary}
\label{cor:seiferter for T}
Let $c$ be a hyperbolic seiferter for $(T_{p, q}, m)$,
where $|p| > q > 2$ and $m \ne pq,\  pq\pm 1$. 
Then, 
\begin{enumerate}
\item
$c$ is $m$--equivalent to a basic seiferter for $T_{p, q}$
$($resp.\ a regular fiber of $S^3 -N(T_{p, q})$$)$ 
if $c$ is an exceptional fiber $($resp.\ a regular fiber$)$
in some Seifert fibration of $T_{p, q}(m)$. 
\item
$c$ cannot be obtained from a basic seiferter or 
a regular fiber of $S^3 -N(T_{p, q})$ by a single $m$--move. 
\end{enumerate}
\end{corollary}

\noindent
\textit{Proof of Corollary~\ref{cor:seiferter for T}}.
It follows from 
the assumption $|p| > q > 2$ and $m \ne pq,\ pq\pm 1$
that $T_{p, q}(m)$ is not a connected sum of
lens spaces, a lens space, or a prism manifold.
Hence, (1) follows from Proposition~\ref{seiferter for T},
and (2) follows from Theorems~\ref{restriction1} and \ref{restriction2}. 
\QED{Corollary~\ref{cor:seiferter for T}}

\bigskip

\bibliographystyle{amsplain}

\begin{thebibliography}{99}

\bibitem{Asano}
K. Asano; 
Homeomorphisms of prism manifolds, 
Yokohama Math.\ J.\ \textbf{26} (1978), 19--25. 

\bibitem{BP}
R. Benedetti and C. Petronio; 
Lectures on hyperbolic geometry, 
Universitext, Springer-Verlag, 1992. 

\bibitem{Berge2}
J. Berge; 
Some knots with surgeries yielding lens spaces, 
Unpublished manuscript. 

\bibitem{BoileauPorti}
M. Boileau and J. Porti; 
Geometrization of 3-orbifolds of cyclic type,  
Ast\'erisque \textbf{272} (2001), 208pp. 

\bibitem{Bonahon}
F. Bonahon and J.-P. Otal; 
Scindements de {H}eegaard des espaces lenticulaires, 
Ann.\ Sci.\ Ecole Norm.\ Sup.\ \textbf{16} (1983), 451--466. 

\bibitem{DEMM}
A. Deruelle, M. Eudave-Mu\~noz, K. Miyazaki and K. Motegi; 
Networking Seifert Surgeries on Knots IV: Seiferters and branched coverings, 
Contemp.\ Math.\ Amer.\ Math.\ Soc.\ \textbf{597} (2013), 235--262. 

\bibitem{DMM2}
A. Deruelle, K. Miyazaki and K. Motegi; 
Networking Seifert Surgeries on Knots II:  Berge's lens surgeries, 
Topology Appl.\ \textbf{156} (2009), 1083--1113. 

\bibitem{DMM1}
A. Deruelle, K. Miyazaki and K. Motegi; 
Networking Seifert Surgeries on Knots, 
Mem.\ Amer.\ Math.\ Soc.\ \textbf{217} (2012),\ no. 1021, viii+130. 

\bibitem{DMMtrefoil}
A. Deruelle, K. Miyazaki and K. Motegi; 
Neighbors of Seifert surgeries on a trefoil knot in the Seifert Surgery Network, 
preprint. 

\bibitem{EM1}
M. Eudave-Mu\~noz;
Non-hyperbolic manifolds obtained by Dehn surgery on a hyperbolic knot,
In: Studies in Advanced Mathematics vol. 2, part~1,
(ed.\ W. Kazez), 1997, Amer.\ Math.\ Soc.\ and International Press,
pp. 35--61.

\bibitem{EM2}
M. Eudave-Mu\~noz; 
On hyperbolic knots with Seifert fibered Dehn surgeries, 
Topology Appl.\ \textbf{121} (2002), 
119--141. 

\bibitem{FS}
R. Fintushel and R. J. Stern; 
Constructing lens spaces by surgery on knots, 
Math.\ Z.\ \textbf{175} (1980), 33--51. 

\bibitem{GW}
C. McA. Gordon and Y.-Q. Wu; 
Annular Dehn fillings, 
Comment.\ Math.\ Helv.\ \textbf{75} (2000), 430--456. 

\bibitem{Greene}
J. E. Greene; 
$L$--space surgeries, genus bound, and the cabling conjecture, 
2012, preprint. 

\bibitem{Hat2}
A. E. Hatcher; 
Notes on basic $3$-manifold topology, 
freely available at \texttt{http://www.math.cornell.edu/hatcher}. 

\bibitem{HiraS}
M. Hirasawa and K. Shimokawa; 
Dehn surgeries on strongly invertible knots which yield lens spaces, 
Proc.\ Amer.\ Math.\ Soc.\ \textbf{128} (2000), 3445--3451. 

\bibitem{HR}
C. Hodgson and J. H. Rubinstein; 
Involutions and isotopies of lens spaces, 
Knot theory and manifolds (Vancouver, B.C., 1983), 
Lect.\ Notes in Math.,\ vol. 1144, Springer-Verlag, 1985, pp.60--96. 

\bibitem{IM}
K. Ishihara and K. Motegi; 
Band sum operations yielding trivial knots, 
Bol.\ Soc.\ Mat.\ Mexicana \textbf{15} (2009), 103--108. 

\bibitem{J}
W. Jaco;
Lectures on three manifold topology, 
CBMS Regional Conference Series in Math., vol. 43,
Amer.\ Math.\ Soc., 1980.

\bibitem{JWW}
B. Jiang and S. Wang and Y.-Q. Wu; 
Homeomorphisms of $3$-manifolds and the realization of Nielsen number, 
Comm.\ Anal.\ Geom.\ \textbf{9} (2001), 825--877. 

\bibitem{MM1}
K. Miyazaki and K. Motegi; 
Seifert fibred manifolds and Dehn surgery, 
Topology\ \textbf{36} (1997), 579--603. 

\bibitem{NR}
W. D. Neumann and F. Raymond; 
Seifert manifolds, plumbing, $\mu$-invariant and orientation reversing maps, 
Algebraic and geometric topology (Proc. Sympos., Univ.  California, Santa Barbara, Calif., 1977), 
Lect.\ Notes in Math.,\ vol. 664, Springer-Verlag, 1978, pp.161--196. 

\bibitem{Oertel}
U. Oertel; 
Closed incompressible surfaces in complements of star links, 
Pacific J.\ Math.\ \textbf{111} (1984), 209--230. 

\bibitem{Orlik}
P. Orlik; 
Seifert manifolds, 
Lect.\ Notes in Math.,\ vol. 291, Springer-Verlag, 1972. 

\bibitem{PetPorti}
C. Petronio and J. Porti;
Negatively oriented ideal triangulations and a proof of Thurston's hyperbolic Dehn filling theorem, 
Expo.\ Math.\ \textbf{18} (2000), 1--35. 

\bibitem{Rubinstein}
J.H. Rubinstein; 
On $3$--manifolds that have finite fundamental groups and contain Klein bottles, 
Trans.\ Amer.\ Math.\ Soc. \textbf{251} 1979, 129--137. 

\bibitem{Sch-Tho}
M. Scharlemann and A. Thompson; 
Link genus and the Conway moves, 
Comment.\ Math.\ Helv.\ \textbf{64} (1989), 527--535. 

\bibitem{T1}
W. P. Thurston; 
The geometry and topology of $3$-manifolds, 
Lecture notes, Princeton University, 1979. 

\bibitem{T2}
W. P. Thurston; 
Three dimensional manifolds, {K}leinian groups and hyperbolic geometry, 
Bull.\ Amer.\ Math.\ Soc.\ \textbf{6} (1982), 357--381. 

\end{thebibliography}

\bigskip

\noindent
Institute of Natural Sciences \par
\noindent
Nihon University \par
\noindent
3-25-40 Sakurajosui, Setagaya-ku \par
\noindent
Tokyo 156-8550 \par
\noindent
Japan \par
\noindent
e-mail: aderuelle@math.chs.nihon-u.ac.jp

\bigskip
\noindent
Faculty of Engineering \par
\noindent
Tokyo Denki University \par
\noindent
5 Senju-Asahi-cho, Adachi-ku \par
\noindent
Tokyo 120--8551 \par
\noindent
Japan \par
\noindent
e-mail: miyazaki@cck.dendai.ac.jp

\bigskip
\noindent
Department of Mathematics \par
\noindent
Nihon University \par
\noindent
3-25-40 Sakurajosui, Setagaya-ku \par
\noindent
Tokyo 156-8550 \par
\noindent
Japan \par
\noindent
e-mail: motegi@math.chs.nihon-u.ac.jp

\end{document}